\documentclass[12pt,leqno]{article}
\usepackage{amsmath,amsfonts,amssymb}
\newtheorem{lem}{Lemma}

\newtheorem{rem}{Remark}
\newtheorem{theo}[lem]{Theorem}

\newtheorem{prop}[lem]{Proposition}
\newcommand{\N}{\mathbb{N}}
\newcommand{\R}{\mathbb{R}}

\newcommand{\E}{\mathbb{E}}
\newcommand{\C}{\mathbb{C}}

\newcommand{\p}{\mathbb{P}}

\newcommand{\pt}{\!\!\!\bf .}

\newcommand{\ind}{\mbox{1}\kern-.21em \mbox{I}}
\newcommand{\wh}{\widehat}
\newcommand{\vp}{\varphi}
\newcommand{\cN}{\mathcal{N}}
\newcommand{\cD}{\mathcal{D}}

\newcommand{\cL}{\mathcal{L}}
\newcommand{\cH}{\mathcal{H}}
\newcommand{\cK}{\mathcal{K}}

\newcommand{\cR}{\mathcal{R}}
\newcommand{\cZ}{\mathcal{Z}}
\font\calcal=cmsy10 scaled\magstep1
\def\build#1_#2^#3{\mathrel{\mathop{\kern 0pt#1}\limits_{#2}^{#3}}}
\def\liml{\build{\longrightarrow}_{}^{{\mbox{\calcal L}}}}

\def\videbox{\mathbin{\vbox{\hrule\hbox{\vrule height1ex \kern.5em\vrule height1ex}\hrule}}}

\def\CO#1{ {\cal O} \left( \frac{1}{T^{#1}} \right)}
\def\esp#1#2{{\mathbb E}_{#1}\left[{#2}\right]}

\def\loga#1{\log \left[ {#1} \right]}
\def\expo#1{\exp \left[ {#1} \right]}

\def\demend{\hfill $\videbox$\\}

\begin{document}
\title{Sharp large deviations for the fractional Ornstein-Uhlenbeck process}
\author{Bernard Bercu
\thanks{Universit\'e Bordeaux 1, Institut de Math\'ematiques de Bordeaux,
UMR C5251, 351 cours de la lib\'eration, 33405 Talence cedex, France. e-mail:
Bernard.Bercu@math.u-bordeaux1.fr}
\and Laure Coutin
\thanks{Universit\'e Paris 5, Centre Universitaire des Saints P\`eres,
UMR C8145,
45 rue des Saints P\`eres, 75270 Paris cedex 06, France. e-mail:
Laure.Coutin@math-info.univ-paris5.fr}
\and Nicolas Savy
\thanks{Universit\'e Paul Sabatier, Institut de Math\'ematiques de Toulouse,
UMR C5583, 31062 Toulouse Cedex 09, France. e-mail:
nicolas.savy@math.univ-toulouse.fr}}
\maketitle
\ \vspace{-1ex} \\
\[\hbox{{\bf Abstract}}\]
We investigate the sharp large deviation properties of the energy and the maximum likelihood estimator
for the Ornstein-Uhlenbeck process driven by a fractional Brownian motion with Hurst index greater than one half.
\vspace{2ex} \newline
{\bf A.M.S. Classification:} 60F10, 60G15, 60J65
\newline
{\bf Key words:} Large deviations, Ornstein-Uhlenbeck process, Fractional brownian motion

\section{Introduction.}
\renewcommand{\theequation}{\thesection.\arabic{equation}}
\setcounter{section}{1}
\setcounter{equation}{0}

Since the pioneer works of Kolmogorov, Hurst and Mandelbrot, a wide range of literature is available on the statistical properties of Fractional Brownian Motion (FBM). On the other hand, one can realize that its large deviations properties were not deeply investigated. Our purpose is to establish sharp large deviations for functionals associated with the Ornstein-Uhlenbeck process driven by a fractional Brownian motion
\begin{eqnarray}
\label{ouf}
dX_t=\theta X_t dt +dW_t^H
\end{eqnarray}
where the initial state $X_0=0$ and the drift parameter $\theta$ is strictly negative. The process $(W^H_t)$ is a FBM  with Hurst parameter $0<H<1$ which means that $(W^H_t)$ is a Gaussian process with continuous paths such that $W_0^H=0$, $\E[W_t^H]=0$ and
\[
\E[W_t^HW_s^H]=\frac{1}{2}\Bigl(t^{2H}+s^{2H} -|t-s|^{2H}\Bigr).
\]
The weighting function $w$ defined, for all $0<s<t$, by
$w(t,s)=w_H^{-1}s^{-H+1/2}(t-s)^{-H+1/2}$
where $w_H$ is a normalizing positive constant, plays a fundamental role for stochastic calculus
associated with $(W^H_t)$. In particular, for $t>0$, let
\begin{eqnarray} \label{martin}
M_t=\int_0^tw(t,s)\,dW_s^H.
\end{eqnarray}
It was proven by Norros {\it et al} \cite{NVV} page 578 that $(M_t)$ is a Gaussian martingale with quadratic variation
$<\!M\!>_t=\lambda_H^{-1}t^{2-2H}$
where
\[
\lambda_H=\frac{8H(1-H)\Gamma(1-2H)\Gamma(H+1/2)}{\Gamma(1/2-H)}
\]
and $\Gamma$ stands for the classical gamma function. In addition, it is also more convenient to study the behavior of
\begin{eqnarray} \label{defY}
Y_t=\int_0^tw(t,s)\,dX_s =\theta\int_0^tw(t,s)X_s\,ds+M_t.
\end{eqnarray}
It was shown by Kleptsyna and Le Breton \cite{KLB1} that whenever, $H>1/2$, relation (\ref{defY}) can be rewritten as
\begin{eqnarray} \label{devY}
Y_t=\theta\int_0^tQ_s\,d\!\!<\!M\!>_s+M_t
\end{eqnarray}
where the process $(Q_t)$ satisfies for all $t>0$
\[
Q_t=\frac{l_H}{2}\left(t^{2H-1}Y_t+\int_0^ts^{2H-1}\,dY_s\right)
\]
with $l_H=\lambda_H/(2(1-H))$. Consequently, we shall assume in all the sequel that $H>1/2$. It follows from \eqref{devY} that the score function, which is the derivative of the log-likelihood function from observations over the interval $[0,T]$, is given by
\[
\Sigma_T(\theta)=\int_0^TQ_t\,dY_t-\theta\int_0^TQ_t^2\,d\!\!<\!M\!>_t.
\]
Via a similar approach as the one of \cite{BR} for the Ornstein-Uhlenbeck process, we shall investigate the large deviation properties for random variables associated with (\ref{ouf}) such as the energy
\begin{eqnarray} \label{defST}
S_T=\int_0^TQ_t^2\,d\!\!<\!M\!>_t
\end{eqnarray}
as well as the maximum likelihood estimator of $\theta$, solution of $\Sigma_T(\theta)=0$, explicitly given by
\begin{eqnarray} \label{defteaT}
\wh{\theta}_T=\frac{\int_0^TQ_t\,dY_t}{\int_0^TQ_t^2\,d\!\!<\!M\!>_t}.
\end{eqnarray}
We also wish to mention the recent work of Bishwal \cite{B08} concerning the large deviation properties of the log-likelihood ratio
$$
(\theta_0 - \theta_1) \int_0^T Q_t dY_t - \frac{(\theta^2_0 - \theta^2_1)}{2} \int_0^T Q^2_t d\!\!<\!M\!>_t.
$$

\noindent As usual, we shall say that a family of real random variables $(Z_T)$ satisfies a Large Deviation Principle (LDP) with rate function $I$, if $I$ is a lower semicontinuous function from $\R$ to $[0,+\infty]$ such that, for any closed set $F \subset \R$,
\begin{eqnarray*} \label{upperb}
\limsup_{T\rightarrow\infty}\frac{1}{T}\log \p(Z_T \in F) \leq -\inf_{x\in F}I(x),
\end{eqnarray*}
while for any open set $G \subset \R$,
\begin{eqnarray*} \label{lowerb}
-\inf_{x\in G}I(x) \leq \liminf_{T\rightarrow\infty}\frac{1}{T}\log \p(Z_T \in G).
\end{eqnarray*}
Moreover, $I$ is a good rate function if its level sets are compact subsets of $\R$. We refer the reader to the excellent book by Dembo and Zeitouni \cite{DZ} on the theory of large deviations, see also \cite{BGR}, \cite{GRZ}. A classical tool for proving an LDP for $S_T$ and $\wh{\theta}_T$ is the normalized cumulant generating function
\begin{eqnarray} \label{defLT}
\cL_T(a,b)=\frac{1}{T} \log \E[\exp(\cZ_T(a,b))]
\end{eqnarray}
where, for any $(a,b)\in\R$,
\begin{eqnarray} \label{defZT}
\cZ_T(a,b)=a\int_0^TQ_t\,dY_t+b\int_0^TQ_t^2\,d\!\!<\!M\!>_t.
\end{eqnarray}
The random variable $\cZ_T(a,b)$ allows us an unified presentation of our results. In fact, in order to establish an LDP for $S_T$ and $\wh{\theta}_T$, it is enough to prove an LDP for $\cZ_T(0,a)$ and $\cZ_T(a,-ca)$, respectively. The following lemma provides an asymptotic expansion for $\cL_T$. It enlightens the role of the limit $\cL$ of $\cL_T$ for the LDP, as well as the first order terms $\cH$ and $\cK_T$ for the sharp LDP. One can observe that it is the keystone of all our results.
\begin{lem}$\pt$ \label{mainlem}
Let $\Delta_H$ be the effective domain of the limit $\cL$ of $\cL_T$
\[
\Delta_H=\Bigl\{(a,b) \in \R^2 \ /\; \theta^2 -2b>0 \ \mbox{and} \
\sqrt{\theta^2 -2b}>\max(a+\theta;-\delta_H(a+\theta))  \Bigr\}
\]
where $\delta_H=(1-\sin(\pi H))/(1+sin(\pi H))$. Then, for any $(a,b)$ in the interior of $\Delta_H$, if $\vp(b)=\sqrt{\theta^2 -2b}$, $\tau(a,b)=\vp(b)-(a+\theta)$ and $r_T(b)=r_H(\vp(b)T/2)\exp(-T\vp(b))-1$, we have the decomposition
\begin{eqnarray} \label{maindeco}
\cL_T(a,b) = \cL(a,b) +\frac{1}{T} \cH(a,b) + \frac{1}{T} \cK_T(a,b) + \frac{1}{T} \cR_T(a,b)
\end{eqnarray}
where
\begin{eqnarray}
\label{defL}
\cL(a,b)&=&-\frac{1}{2}(a+\theta +\vp(b)),\\
\label{defH}
\cH(a,b)&=&-\frac{1}{2}\log\left(\frac{\tau(a,b)}{2\vp(b)}\right),\\
\label{defK}
\cK_T(a,b)&=&-\frac{1}{2}\log\left(1+\frac{(2\vp(b)-\tau(a,b))}{2\vp(b)}r_T(b)\right).
\end{eqnarray}
Here, $I_H$ is the modified Bessel function of the first kind \cite{Leb} and the function $r_H$ is defined for all $z \in \C$ with $|\arg z|<\pi$ by
\begin{eqnarray*}
r_H(z)=\frac{\pi z}{\sin(\pi H)} \Bigl(I_H(z)I_{1-H}(z)+I_{-H}(z)I_{H-1}(z)\Bigr).
\end{eqnarray*}
Finally, the remainder is
\begin{eqnarray} \label{defR}
\cR_T(a,b)=-\frac{1}{2} \log \left(1 + \frac{(2\vp(b)-\tau(a,b))^2} {\tau(a,b)(2\vp(b)+r_T(b)(2\vp(b)-\tau(a,b)))} e^{-2T\vp(b)}\! \right)\!\!.
\end{eqnarray}
\end{lem}

\noindent{\bf Proof.} The proof is given in Appendix A.\demend

\begin{rem}$\pt$
By use of the duplication formula for the gamma function \cite{Leb}, one can realize that if $H=1/2$,
$r_H(z)=e^{2z}+e^{-2z}$ which immediately leads to $r_T(b)=e^{-2T\vp(b)}$. Consequently, in that particular case, $\cK_T(a,b)$ as well as $\cR_T(a,b)$ go exponentially fast to zero and we find again Lemma 2.1 of \cite{BR} which is the keystone for all results in \cite{BR}.
\end{rem}

\section{The energy.}
\renewcommand{\theequation}{\thesection.\arabic{equation}}
\setcounter{section}{2}
\setcounter{equation}{0}
First of all, we shall focus our attention on the energy $(S_T)$. One can observe that the strong law of large numbers as well as the Central Limit Theorem (CLT) for the sequence $(S_T)$ were not previously established in the literature.
\begin{prop}$\pt$ \label{lemenergy}
We have
\begin{equation} \label{energyslln}
\lim_{T\rightarrow \infty}\frac{S_T}{T}= -\frac{1}{2\theta}\hspace{1.5cm}\mbox{a.s.}
\end{equation}
Moreover, we also have the CLT
\begin{equation} \label{energyclt}
\frac{1}{\sqrt{T}}\left(S_T+\frac{T}{2\theta}\right)\liml \cN\left(0,-\frac{1}{2\theta^3}\right).
\end{equation}
\end{prop}

\noindent{\bf Proof of Proposition \ref{lemenergy}.}
The almost sure convergence (\ref{energyslln}) clearly follows from Theorem \ref{theoldpener} below as
the sequence $(S_T/T)$ satisfies an LDP with good rate function $I$ given by (\ref{ratenergy}). It is not hard to see that $I(c)=0$ if and only if $c=-1/2\theta$. In order to prove the CLT given by (\ref{energyclt}), denote
$$
V_T=\frac{1}{\sqrt{T}}\left(S_T+\frac{T}{2\theta}\right).
$$
For all $a \in \R$, let $\Lambda_T(a)=\E[\exp(a V_T)]$. We clearly have
\begin{equation*}
\Lambda_T(a)=\exp\left(\frac{a\sqrt{T}}{2\theta}\right)\E\left[\exp\Bigl(\frac{aS_T}{\sqrt{T}}\Bigr)\right].
\end{equation*}
Hence, we deduce from the decomposition (\ref{maindeco}) with $a=0$ and $b=a$ that
\begin{equation} \label{proofCLTenergy1}
\Lambda_T(a)=\exp\!\left(\!\frac{a\sqrt{T}}{2\theta} + T \cL\Bigl(0,\frac{a}{\sqrt{T}}\Bigr) +
\cH\Bigl(0,\frac{a}{\sqrt{T}}\Bigr) + \cK_T\Bigl(0,\frac{a}{\sqrt{T}}\Bigr) + \cR_T\Bigl(0,\frac{a}{\sqrt{T}}\Bigr)\!\right)\!.
\end{equation}
On the one hand,
$$
\cL\Bigl(0,\frac{a}{\sqrt{T}}\Bigr)=-\frac{1}{2}(\theta +\vp_T)
$$
where
$$
\vp_T=\sqrt{\theta^2 -\frac{2a}{\sqrt{T}}}=-\theta\sqrt{1 -\frac{2a}{\theta^2\sqrt{T}}}.
$$
Consequently, as
$$
\vp_T=-\theta + \frac{a}{\theta \sqrt{T}} + \frac{a^2}{2\theta^3 T} + o\left(\frac{1}{T}\right)
$$
we obtain that
\begin{equation} \label{proofCLTenergy2}
\lim_{T\rightarrow \infty}\frac{a\sqrt{T}}{2\theta} + T \cL\Bigl(0,\frac{a}{\sqrt{T}}\Bigr) = -\frac{a^2}{4\theta^3}.
\end{equation}
On the other hand, as $(\vp_T)$ converges to $-\theta$, one can check that
\begin{equation} \label{proofCLTenergy3}
\lim_{T\rightarrow \infty} \left( \cH\Bigl(0,\frac{a}{\sqrt{T}}\Bigr)+ \cK_T\Bigl(0,\frac{a}{\sqrt{T}}\Bigr) + \cR_T\Bigl(0,\frac{a}{\sqrt{T}}\Bigr) \right)=0.
\end{equation}
Therefore, we infer from (\ref{proofCLTenergy1}), (\ref{proofCLTenergy2}) and (\ref{proofCLTenergy3}) that
\begin{equation} \label{proofCLTenergy4}
\lim_{T\rightarrow \infty} \Lambda_T(a)=\exp\left(-\frac{a^2}{4\theta^3}\right).
\end{equation}
Convergence (\ref{proofCLTenergy4}) holds for all $a$ in a neighborhood of the origin which leads to (\ref{energyclt}) and completes the proof of Proposition \ref{lemenergy}.\demend

\noindent
In order to obtain the large deviation properties for $(S_T)$, we shall make use of Lemma \ref{mainlem} with $a=0$ and $b=a$. On the one hand, let
\[
D_H=\Bigl\{a \in \R \ /\; \theta^2 -2a>0 \ \mbox{and} \ \sqrt{\theta^2 -2a}>-\delta_H\theta  \Bigr\}.
\]
It is not hard to see that $D_H=]-\infty,a_H[$ where
\begin{equation} \label{defaHenergy}
a_H=\frac{\theta^2}{2}(1-\delta_H^2).
\end{equation}
Consequently, as $|\delta_H|<1$, one can observe that the origin always belongs to the interior of $D_H$. On the other hand, for all $a \in D_H$, let $\vp(a)=\sqrt{\theta^2 -2a}$,
\begin{eqnarray}
\label{defLS}
L(a)&=&\cL(0,a)=-\frac{1}{2}(\theta+\sqrt{\theta^2 -2a}),\\
\label{defHS}
H(a)&=&\cH(0,a)=-\frac{1}{2}\log\left(\frac{\vp(a)-\theta}{2\vp(a)}\right),\\
\label{defKTS}
K_T(a)&=&\cK_T(0,a)= -\frac{1}{2}\log\left(1+\frac{(\vp(a)+\theta)}{2\vp(a)}r_T(a)\right).
\end{eqnarray}
The main difficulty comparing to \cite{BR} is that the function $L$ is not steep. Actually,
\[
L^\prime(a)=\frac{1}{2\sqrt{\theta^2 -2a}}
\]
which implies that $L^\prime(a_H)=-1/(2\theta\delta_H)$. Moreover, for all $c>0$, $L^\prime(a)=c$ if and only if  $a=a_c$ with $a_c=(4\theta^2c^2 -1)/(8c^2)$. Hence, $a_c<a_H$ whenever $0<c<-1/(2\theta\delta_H)$. Denote by $I$ the Fenchel-Legendre transform of the function $L$. Our first result on the large deviation properties of $(S_T/T)$ is as follows.

\begin{theo}$\pt$ \label{theoldpener}
The sequence $(S_T/T)$ satisfies an LDP with good rate function
\begin{eqnarray} \label{ratenergy}
I(c)= \left \{ \begin{array}{ll}
    {\displaystyle \frac{(2\theta c +1)^2}{8c}\ \mbox{ if }\ 0<c\leq -\frac{1}{2\theta \delta_H},} \vspace{2ex}\\
    {\displaystyle \frac{c\theta^2}{2}(1-\delta_H^2)+\frac{\theta}{2}(1-\delta_H)\ \mbox{ if }\ c\geq -\frac{1}{2\theta \delta_H},} \vspace{2ex}\\
    {\displaystyle +\infty\ \mbox{ otherwise. }}
   \end{array}  \right.
\end{eqnarray}
\end{theo}

\begin{rem}$\pt$
In the particular case $H=1/2$, then $\delta_H=0$ and the LDP for $(S_T/T)$ is exactly the one established by Bryc and Dembo \cite{BD} for general centered Gaussian processes.
\end{rem}

We are now going to improve this result by a Sharp Large Deviation Principle (SLDP) for $(S_T/T)$ inspired by the well-known Bahadur-Rao Theorem \cite{BaR} on the sample mean.


\begin{theo}$\pt$  \label{theosldpener}
The sequence $(S_T/T)$ satisfies a SLDP associated with $L$, $H$ and $K_T$ given by \eqref{defLS}, \eqref{defHS}, and \eqref{defKTS}, respectively.
\newline
a) For all $- 1/(2\theta)< c <- 1/(2\theta \delta_H)$, there exists a sequence $(b_{c,k}^H)$
such that, for any $p>0$ and $T$ large enough,
\begin{equation} \label{mainsldener1}
\p(S_T \geq cT) =  \frac{\exp(-TI(c) +  J(c) + K_H(c))}{a_c \sigma_c \sqrt{2 \pi T}}
\left[1 + \sum_{k=1}^{p}\frac{b_{c,k}^H}{T^k} + {\cal O}\Bigl(\frac{1}{T^{p+1}}\Bigr)\right]
\end{equation}
while, for $0<c<- 1/(2\theta)$,
\begin{equation} \label{mainsldener2}
\p(S_T \leq cT) =  - \frac{\exp(-TI(c) +  J(c) + K_H(c))}{a_c \sigma_c \sqrt{2 \pi T}}
\left[1 + \sum_{k=1}^{p}\frac{b_{c,k}^H}{T^k} + {\cal O}\Bigl(\frac{1}{T^{p+1}}\Bigr)\right]
\end{equation}
where
\begin{equation} \label{defJKener1}
J(c)=-\frac{1}{2} \log \left(\frac{1-2\theta c}{2}\right),
\hspace{1cm}
K_H(c)=-\frac{1}{2} \log \left(\frac{(1+ \sin(\pi H))(1+2\theta c\delta_H)}{ 2 \sin(\pi H) }\right),
\end{equation}
with
\begin{equation} \label{defacsgener1}
a_c=\frac{4\theta^2c^2 -1}{8c^2}
\hspace{1.5cm}
\text{and}
\hspace{1.5cm}
\sigma_c^2=4c^3.
\end{equation}
Moreover, the coefficients $b_{c,1}^H, b_{c,2}^H, \ldots,b_{c,p}^H$ may be explicitly calculated as functions of the derivatives of $L$ and $H$ evaluated at point $a_c$. They also depend on the Taylor expansion of $K_T$ and its derivatives at $a_c$.
\newline
b) For all $c  > - 1/(2\theta \delta_H)$, there exists a sequence $(d_{c,k}^H)$ such that, for any $p>0$ and $T$ large enough
\begin{equation} \label{mainsldener3}
\p(S_T \geq cT) =  \frac{\exp(-TI(c) +  P_H(c) + Q_H(c))}{a_H \sigma_H \sqrt{2 \pi T}}
\left[1 + \sum_{k=1}^{p}\frac{d_{c,k}^H}{T^k} + {\cal O}\Bigl(\frac{1}{T^{p+1}}\Bigr)\right]
\end{equation}
\begin{equation} \label{defPQener3}
P_H(c)=-\frac{1}{2} \log \left(\frac{-(1+2\theta c\delta_H)}{4\delta_H \sin(\pi H)}\right),
\hspace{1cm}
Q_H(c)=\frac{(2H-1)^2\sin(\pi H)(1+2\theta c\delta_H)}{2(1-(\sin(\pi H))^2)}
\end{equation}
with
\begin{equation}
\label{defaHsHgener3}
a_H=\frac{\theta^2(1-\delta_H^2)}{2}
\hspace{1.5cm}
\text{and}
\hspace{1.5cm}
\sigma_H^2=-\frac{1}{2\theta^3\delta_H^3}.
\end{equation}
Moreover, the coefficients $d_{c,1}^H, d_{c,2}^H, \ldots,d_{c,p}^H$ may be explicitly calculated
as functions of the derivatives of $L$ and $H$ evaluated at point $a_H$. They also depend on the Taylor expansion of $K_T$ and its derivatives at $a_H$.
\end{theo}

\begin{rem}$\pt$
For example, the first coefficient $b_{c,1}^H$ is given by
\begin{equation*}
b_{c,1}^H=\frac{1}{\sigma_c^2}\left(
\frac{s_1}{a_c} -\frac{s_1^2}{2} - \frac{s_2}{2}
-\frac{s_3}{2a_c\sigma_c^2}
+\frac{s_1\ell_3}{2\sigma_c^2}-\frac{5\ell_3^2}{24\sigma_c^4}
+\frac{\ell_4}{8\sigma_c^2} -\frac{1}{a_c^2}
\right) + k_{c,1}^H
\end{equation*}
where
$\ell_q=L^{(q)}(a_c)$, $h_q=H^{(q)}(a_c)$, $s_q=h_q+k_q$, $p_H=(1- \sin(\pi H))/sin(\pi H)$,  with
\begin{eqnarray*}
k_1&=&\lim_{T\rightarrow \infty}K_T^{\prime}(a_c)=\frac{-4\theta p_H c^3}{2+p_H(1+2\theta c)},\\
k_2&=&\lim_{T\rightarrow \infty}K_T^{\prime\prime}(a_c)=\frac{16\theta p_H c^5(6+p_H(3+2\theta c))}{(2+p_H(1+2\theta c))^2},\\
k_{c,1}^H&=&\lim_{T\rightarrow \infty}T(K_T(a_c)-K_H(c))=\frac{c(1+2\theta c)(2H-1)^2}{2\sin(\pi H)(2+p_H(1+2\theta c))}.
\end{eqnarray*}
In addition, one can also observe that $\sigma_c^2=L^{\prime\prime}(a_c)$ and $\sigma_H^2=L^{\prime\prime}(a_H)$.
\end{rem}


\begin{theo} \label{theosldpenereg}
For $c  = - 1/(2\theta \delta_H)$, there exists a sequence $(d_k^H)$
such that, for any $p>0$ and $T$ large enough
\begin{equation} \label{mainsldener4}
\p(S_T \geq cT) = \frac{\exp(-TI(c)+K_H)\Gamma(1/4)}{2 \pi a_H \sigma_H  T^{1/4}}
\left[1 + \sum_{k=1}^{2p} \frac{d_k^H}{(\sqrt{T})^k} + {\cal O}\Bigl(\frac{1}{T^{p}\sqrt{T}}\Bigr)\right]
\end{equation}
where $a_H$ and $\sigma^2_H$ are given by \eqref{defaHsHgener3} and
$$
K_H=\frac{1}{2}\log(\delta_H \sin(\pi H))+\frac{1}{4}\log(-\theta \delta_H).
$$
As before, the coefficients $d_1^H, d_2^H, \ldots, d_{2p}^H$ may be explicitly calculated.
\end{theo}

\noindent{\bf Proof.} The proofs are given in Section \ref{S-Proof}. \demend

\section{The maximum likelihood estimator.}
\renewcommand{\theequation}{\thesection.\arabic{equation}}
\setcounter{section}{3}
\setcounter{equation}{0}

Our purpose is now to establish similar results for the maximum likelihood estimator. 
The almost sure convergence of $\wh{\theta}_T$ towards $\theta$
was already proven by  Kleptsyna and Le Breton \cite{KLB1}, see also
Prakasa Rao \cite{PKS1, PKS2} for related results. We just learn that
the central limit theorem was established via a different approach
by Brouste and Kleptsyna \cite{BK}. 
\begin{prop}$\pt$ \label{propmle}
We have
\begin{equation} \label{mleslln}
\lim_{T\rightarrow \infty}\wh{\theta}_T = \theta\hspace{1.5cm}\mbox{a.s.}
\end{equation}
Moreover, we also have the CLT
\begin{equation} \label{mleclt}
\sqrt{T}\Bigl(\wh{\theta}_T-\theta\Bigr)\liml \cN\Bigl(0,-2\theta\Bigr).
\end{equation}
\end{prop}

\noindent{\bf Proof of Proposition \ref{propmle}.}
The almost sure convergence (\ref{mleslln}) is given by Proposition 2.2 of Kleptsyna and Le Breton \cite{KLB1}. For all $c\in \R$, denote
$$
V_T(c)=\frac{1}{\sqrt{T}}\int_0^T Q_tdY_t-\left(\frac{c}{\sqrt{T}}+\theta \right)\frac{S_T}{\sqrt{T}}.
$$
One can easily check that
\begin{equation} \label{proofCLTmle1}
\p(\sqrt{T}(\wh{\theta}_T- \theta) \leq c)=\p(V_T(c) \leq 0).
\end{equation}
Consequently, in order to prove the CLT given by (\ref{mleclt}), it is only necessary to establish the asymptotic behavior the sequence $(V_T(c))$ where the threshold $c$ can be seen as a parameter. For all $a,c \in \R$, let $\Lambda_T(a,c)=\E[\exp(a V_T(c))]$.
We clearly have
\begin{equation*}
\Lambda_T(a,c)=\exp\left(T L_T\Bigl(\frac{a}{\sqrt{T}}, c_T\Bigr)\right)
\hspace{0.5cm}\text{where}\hspace{0.5cm}
c_T=-\frac{a}{\sqrt{T}}\left(\frac{c}{\sqrt{T}}+\theta \right).
\end{equation*}
Thus, it follows from the decomposition (\ref{maindeco}) that
\begin{equation} \label{proofCLTmle2}
\Lambda_T(a,c)=\exp\!\left(T \cL\Bigl(\frac{a}{\sqrt{T}}, c_T\Bigr) +
\cH\Bigl(\frac{a}{\sqrt{T}}, c_T\Bigr)+
\cK_T\Bigl(\frac{a}{\sqrt{T}}, c_T\Bigr)+\cR_T\Bigl(\frac{a}{\sqrt{T}}, c_T\Bigr)
\!\right)\!.
\end{equation}
On the one hand,
$$
\cL\Bigl(\frac{a}{\sqrt{T}}, c_T\Bigr)=-\frac{1}{2}\left(\frac{a}{\sqrt{T}}+\theta +\vp_T\right)
$$
where
$$
\vp_T=\sqrt{\theta^2 -2c_T}=-\theta\sqrt{1 -\frac{2c_T}{\theta^2}}.
$$
Hence, as
$$
\vp_T=-\theta - \frac{a}{\sqrt{T}} + \frac{(a^2-2 a c)}{2\theta T} + o\left(\frac{1}{T}\right)
$$
we deduce that
\begin{equation} \label{proofCLTmle3}
\lim_{T\rightarrow \infty} T \cL\Bigl(\frac{a}{\sqrt{T}}, c_T\Bigr) = -\frac{(a^2-2 a c)}{4\theta}.
\end{equation}
On the other hand, as $(\vp_T)$ converges to $-\theta$, it is not hard to see that
\begin{equation} \label{proofCLTmle4}
\lim_{T\rightarrow \infty} \left( \cH\Bigl(\frac{a}{\sqrt{T}}, c_T\Bigr)+
\cK_T\Bigl(\frac{a}{\sqrt{T}}, c_T\Bigr)+\cR_T\Bigl(\frac{a}{\sqrt{T}}, c_T\Bigr) \right) =0.
\end{equation}
The conjunction of (\ref{proofCLTmle2}), (\ref{proofCLTmle3}) and (\ref{proofCLTmle4}) leads to
\begin{equation} \label{proofCLTmle5}
\lim_{T\rightarrow \infty} \Lambda_T(a,c)=\exp\left(-\frac{(a^2-2 a c)}{4\theta}\right).
\end{equation}
This convergence holds for all $a$ in a neighborhood of the origin. Consequently,
\begin{equation} \label{proofCLTmle6}
V_T(c) \liml \cN\left(\frac{c}{2\theta},\frac{-1}{2\theta}\right).
\end{equation}
Denote by $V(c)$ the limiting distribution of $(V_T(c))$. It follows from a standard Gaussian calculation that
\begin{equation} \label{proofCLTmle7}
\p(V(c)\leq 0)= \frac{1}{-4\pi \theta}\int_{-\infty}^c\exp\Bigl(\frac{x^2}{2\theta}\Bigr)\,dx.
\end{equation}
Finally, (\ref{proofCLTmle1}) and (\ref{proofCLTmle7}) imply (\ref{mleclt}) which achieves the proof of Proposition \ref{propmle}. \demend

\noindent In order to establish the large deviation properties of $(\wh{\theta}_T)$, we shall make use of the auxiliary random variable defined for all $c \in \R$ by
\begin{eqnarray} \label{link}
Z_T(c)=\int_0^T Q_tdY_t-c\int_0^T Q_t^2\,d\!\!<\!M\!>_t
\end{eqnarray}
since $\p(\wh{\theta}_T \geq c)=\p(Z_T(c) \geq 0)$. Let
\[
D_H=\Bigl\{a \in \R \ /\; \theta^2 +2ac>0 \ \mbox{and} \ \sqrt{\theta^2 +2ac} > \max(a+\theta;-\delta_H(a+\theta))  \Bigr\}.
\]
After some straightforward calculations, it is not hard to see that
\begin{eqnarray*}
D_H= \left \{ \begin{array}{lll}
    {\displaystyle \  ]a_1^H, a_2^H[ \ \ \mbox{ if }\ c \leq
    \frac{\theta}{2}}, \vspace{1ex}\\
    {\displaystyle  \  ]a_1^H,  a^c[
    \ \ \ \mbox{ if }\ c>\frac{\theta}{2}},
   \end{array}  \right.
\end{eqnarray*}
where $a^c=2(c-\theta)$ and
\begin{eqnarray*}
\left \{ \begin{array}{lll}
    {\displaystyle a_1^H=\frac{c-\theta
    \mu_H-\sqrt{c^2-2\theta c \mu_H+\theta^2\mu_H}}{\mu_H}},\vspace{1ex}\\
    {\displaystyle a_2^H=\frac{c-\theta
    \mu_H+\sqrt{c^2-2\theta c \mu_H+\theta^2\mu_H}}{\mu_H}},
   \end{array}  \right.
\end{eqnarray*}
where $\mu_H = \delta_H^2$. In addition, for all $a\in D_H$, let $\vp(a) = \sqrt{\theta^2 +2ac}$ and
\begin{eqnarray}
\label{defLmle}
L(a)&=&\cL(a,-ca)= -\frac{1}{2}\left(a+\theta +\sqrt{\theta^2 +2ac}\right),\\
\label{defHmle}
H(a)&=&\cH(a,-ca)=-\frac{1}{2}\log\left(\frac{\vp(a)-a-\theta}{2\vp(a)}\right),\\
\label{defKTmle}
K_T(a)&=&\cK_T(a,-ca)= -\frac{1}{2}\log\left(1+\frac{(\vp(a)+a+\theta)}{2\vp(a)}r_T(-ac)\right).
\end{eqnarray}

The function $L$ is not steep as the derivative of $L$ is finite at the boundary of $D_H$. Moreover, $L^{'}(a)=0$ if and only if $a=a_c$ with $a_c=(c^2-\theta^2)/(2c)$ and $a_c \in D_H$ whenever $c<\theta/3$.

\begin{theo} $\pt$ \label{theoldpmle}
The maximum likelihood estimator $(\wh{\theta}_T)$ satisfies an LDP with good rate function
\begin{eqnarray} \label{ratemle}
I(c)=\left \{ \begin{array}{ll}
    {\displaystyle \!-\frac{(c-\theta)^2}{4c}
    \ \ \mbox{ if }\ c<\frac{\theta}{3}}, \vspace{1ex}\\
    {\displaystyle \ \ 2c-\theta \ \ \ \ \mbox{ if }\ c\geq \frac{\theta}{3}.}
   \end{array}  \right.
\end{eqnarray}
\end{theo}

\begin{rem}$\pt$
One can observe that the rate function $I$ is totally free of the parameter $H$. Hence, $(\wh{\theta}_T)$ shares the same LDP than the one established by Florens-Landais and Pham \cite{FlP} for the standard Ornstein-Uhlenbeck process with $H=1/2$.
\end{rem}

\begin{theo}$\pt$ \label{theosldpmle}
The maximum likelihood estimator $(\hat{\theta}_T)$ satisfies an SLDP associated with $L$, $H$ and $K_T$ given by \eqref{defLmle}, \eqref{defHmle}, and \eqref{defKTmle}, respectively.
\newline
a) For all $\theta < c <\theta/3$, there exists a sequence $(b_{c,k}^H)$ such that, for any $p>0$ and $T$ large enough,
\begin{equation} \label{mainsldpmle1}
\p(\hat{\theta}_T \geq c) = \frac{\exp(-TI(c) +  J(c)+K_H(c))}{\sigma_c a_c\sqrt{2 \pi T}}
\left[1 + \sum_{k=1}^{p}\frac{b_{c,k}^H}{T^k} + {\cal O}\Bigl(\frac{1}{T^{p+1}}\Bigr)\right]
\end{equation}
while, for $c<\theta$,
\begin{equation} \label{mainsldpmle2}
\p(\hat{\theta}_T \leq c) = - \frac{\exp(-TI(c) +  J(c)+K_H(c))}{\sigma_c a_c\sqrt{2 \pi T}}
\left[1 + \sum_{k=1}^{p}\frac{b_{c,k}^H}{T^k} + {\cal O}\Bigl(\frac{1}{T^{p+1}}\Bigr)\right]
\end{equation}
where
\begin{equation} \label{defHacmle1}
J(c)=-\frac{1}{2} \log \left(\frac{(c+\theta)(3c-\theta)}{4c^2}\right),
\hspace{1cm}
K_H(c)=-\frac{1}{2} \log \left(1+p_H\frac{(c-\theta)^2}{4c^2}\right)
\end{equation}
with $p_H=(1- \sin(\pi H))/sin(\pi H)$,
\begin{equation} \label{defacsgmle1}
a_c=\frac{c^2-\theta^2}{2c}
\hspace{1.5cm}
\mbox{and}
\hspace{1.5cm}
\sigma_c^2=-\frac{1}{2c}.
\end{equation}
Moreover, the coefficients $b_{c,1}^H, b_{c,2}^H, \ldots,b_{c,p}^H$ may be explicitly calculated as in Theorem \ref{theosldpener}.
\newline
b) For all $c > \theta/3$ with $c \ne 0$, 
there exists a sequence $(d_{c,k}^H)$ such that, for any $p> 0$ and $T$ large enough,
\begin{equation} \label{mainsldpmle3}
\p(\hat{\theta}_T \geq c) =  \frac{\exp(-TI(c) +  P(c))\sqrt{\sin(\pi H)}}{\sigma^c a^c\sqrt{2 \pi T}}
\left[1 + \sum_{k=1}^{p}\frac{d_{c,k}^H}{T^k} + {\cal O}\Bigl(\frac{1}{T^{p+1}}\Bigr)\right]
\end{equation}
where
\begin{equation} \label{defPacmle}
P(c)=-\frac{1}{2} \log \left(\frac{(c-\theta)(3c-\theta)}{4c^2}\right)
\end{equation}
with
\begin{equation}
\label{defacsgmle2}
a^c=2(c-\theta)
\hspace{1.5cm}
\mbox{and}
\hspace{1.5cm}
(\sigma^c)^2=\frac{c^2}{2(2c-\theta)^3}.
\end{equation}
In addition, the coefficients $d_{c,1}^H, d_{c,2}^H, \ldots,d_{c,p}^H$ may be explicitly calculated as in Theorem \ref{theosldpener}.
\newline
c) Finally, for $c=0$, for any $p>0$ and for $T$ large enough,
\begin{equation} \label{mainsldpmle5}
\p(\hat{\theta}_T \geq 0) = 2 \, \frac{\exp(-TI(c))\sqrt{\sin(\pi H)}}
{\sqrt{2 \pi T}\sqrt{-2\theta}}
\left[1 + \sum_{k=1}^{p}\frac{d_{k}^H}{T^k} + {\cal O}\Bigl(\frac{1}{T^{p+1}}\Bigr)\right]
\end{equation} 
As before, the coefficients $d_1^H, d_2^H, \ldots, d_{p}^H$ may be explicitly calculated. 
\end{theo}

\begin{theo}
For $c = \theta/3$, there exists a sequence $(d_k^H)$ such that, for any $p>0$ and $T$ large enough 
\begin{equation} \label{mainsldpmle4}
\p(\hat{\theta}_T \geq \frac{\theta}{3}) = \frac{\exp(-TI(c)) \, \Gamma\left( \frac{1}{4} \right)}{4 \pi T^{1/4} \, a_{\theta}^{3/4} \, \sigma_{\theta}} \, \sqrt{\sin(\pi H)} \left[1 + \sum_{k=1}^{p}\frac{e_k^H}{T^k} + {\cal O}\Bigl(\frac{1}{T^{p+1}}\Bigr)\right] 
\end{equation} 
where $a_{\theta}$ and $\sigma_{\theta}$ are given by 
\begin{equation} 
a_{\theta}=-\frac{4 \theta}{3} 
\hspace{1.5cm} 
\mbox{and} 
\hspace{1.5cm} 
\sigma_{\theta}^2=-\frac{3}{2 \theta}. 
\end{equation} 
As before, the coefficients $e_{1}^H, e_{2}^H, \ldots,e_{p}^H$  may be explicitly calculated.
\end{theo} 

\noindent{\bf Proof.} 
The proofs are left to the reader as they follow essentially the same lines as that for the energy.\demend 

\section{Proofs of the main results.} \label{S-Proof}
\subsection{Proof of Theorem \ref{theosldpener}, first part.}
We first focuse our attention on the sharp LDP for the energy $S_T$ in the easy case $- 1/(2\theta)<c<- 1/(2\theta \delta_H)$. Let $L_T$ be the normalized cumulant generating function of $S_T$. We already saw that $a_c$, given by \eqref{ratenergy}, belongs to the domain $D_H$ whenever $c<-1/(2\theta\delta_H)$. For all $- 1/(2\theta)<c<- 1/(2\theta \delta_H)$, we can split $\p(S_T \geq cT)$ into two terms, $\p(S_T \geq cT) = A_TB_T$ with
\begin{eqnarray}
\label{defAener}
A_T&=&\exp(T(L_T(a_c)-ca_c)), \\
\label{defBener}
B_T&=&\E_{T}\Big[\exp(-a_c(S_T-cT))\ind_{S_T\geq cT}\Big],
\end{eqnarray}
where $\E_T$ stands for the expectation after the usual change of probability
\begin{eqnarray}
\label{newprobstan}
\frac{d\p_T}{d\p} = \exp\Bigl(a_cS_T- TL_T(a_c)\Bigr).
\end{eqnarray}
On the one hand, we can deduce from (\ref{maindeco}) with $L(a)={\cal L}(0,a)$, $H(a)={\cal H}(0,a)$, $K_T(a)=\cK_T(0,a)$, and $R_T(a)=\cR_T(0,a)$ together with (\ref{ratenergy}) and (\ref{defJKener1}) that
\begin{eqnarray}
A_T & = &\exp\Bigl(T(L(a_c) - c a_c) + H(a_c)+K_T(a_c)+R_T(a_c)\Bigr),\nonumber\\
\label{devAener}
A_T & = &\exp\Bigl(-TI(c)+J(c)+K_T(a_c)+R_T(a_c)\Bigr).
\end{eqnarray}
Consequently, we infer from (\ref{defKTS}) and (\ref{defJKener1}) that for any $p>0$ and $T$ large enough
\begin{eqnarray} \label{devKener}
K_T(a_c)=K_H(c) + \sum_{k=1}^{p}\frac{\gamma_{k}}{T^k} + {\cal O}\Bigl(\frac{1}{T^{p+1}}\Bigr)
\end{eqnarray}
where the coefficients  $(\gamma_{k})$, which also depend on $H$, may be explicitly calculated. In addition, it is not hard to see from  (\ref{defR}) that the remainder $R_T(a_c)={\cal O}(\exp(-T/c))$.
Therefore, we deduce from (\ref{devAener}) and (\ref{devKener}) that for any $p>0$ and $T$ large enough \begin{eqnarray} \label{devfinAener}
A_T=\exp(-TI(c) + J(c)+K_H(c)) \left[1 + \sum_{k=1}^{p}\frac{\alpha_{k}}{T^k} + {\cal O}\Bigl(\frac{1}{T^{p+1}}\Bigr)\right]
\end{eqnarray}
where, as before, the coefficients $(\alpha_{k})$ may be explicitly calculated. For example, $$
\alpha_1= \frac{-2c(1+2\theta c)r_{1}^H}{\sin(\pi H)(2+(1+2\theta c)p_H)}.
$$
It now remains to give the expansion for $B_T$ which can be rewritten as
\begin{equation} \label{devBener}
B_T = \E_T \Big[\exp(-a_c \sigma_c \sqrt{T} U_T) \ind_{U_T\geq 0}\Big]
\end{equation}
where
$$
U_T = \frac{S_T -cT}{\sigma_c\sqrt{T}}.
$$

\begin{lem}$\pt$
\label{lemBenereasy}
For all $- 1/(2\theta)<c<- 1/(2\theta \delta_H)$, the distribution of $U_T$ under $\p_T$ converges, as $T$ goes to infinity, to the $\cN(0,1)$ distribution. Moreover, there exists a sequence $(\beta_k)$ such that, for any $p>0$ and $T$ large enough,
\begin{equation} \label{devfinBener}
B_T= \frac{\beta_0}{\sqrt{T}}\left[ 1+\sum_{k=1}^p \frac{\beta_k}{T^k} + {\mathcal O}\Bigl( \frac{1}{T^{p+1}}\Bigr)\right].
\end{equation}
The sequence $(\beta_k)$ only depends on the derivatives of $L$ and $H$ evaluated at point $a_c$. They also depend on the Taylor expansion of $K_T$ and its derivatives at $a_c$. For example,
\begin{equation*}
\beta_0=\frac{1}{a_c\sigma_c\sqrt{2\pi}},
\end{equation*}
\vspace{-1ex}
\begin{equation*}
\beta_1=\frac{1}{\sigma_c^2}\left(
\frac{s_1}{a_c} -\frac{s_1^2}{2} - \frac{s_2}{2}
-\frac{s_3}{2a_c\sigma_c^2}
+\frac{s_1\ell_3}{2\sigma_c^2}-\frac{5\ell_3^2}{24\sigma_c^4}
+\frac{\ell_4}{8\sigma_c^2} -\frac{1}{a_c^2}
\right),
\end{equation*}
where
$\ell_q=L^{(q)}(a_c)$, $h_q=H^{(q)}(a_c)$, and $s_q=h_q+k_q$.
\end{lem}

\noindent
{\bf Proof.} The proof of Lemma \ref{lemBenereasy} is given in Appendix B.\demend

\noindent
{\bf Proof of Theorem \ref{theosldpener}, first part.} The expansions \eqref{mainsldener1} and (\eqref{mainsldener2} immediately follow from the conjunction of \eqref{devfinAener} and \eqref{devfinBener}.
\demend

\subsection{Proof of Theorem \ref{theosldpener}, second part.}

We are now going to establish the sharp LDP for the energy $S_T$ in the more complicated case $c>- 1/(2\theta \delta_H)$. We have already seen that the point $a_c$, given by \eqref{ratenergy}, belongs
to the domain $D_H=]-\infty, a_H[$ whenever $c<-1/(2\theta\delta_H)$. Consequently, for $c>- 1/(2\theta \delta_H)$, it is necessary to make use of a slight modification of the strategy of time varying change of probability proposed by Bryc and Dembo \cite{BD}.\\

Now, we show a modification of the decomposition \eqref{maindeco} which  allows us to replace $r_T(a)$ by $p_H$ in $K_T$, putting the difference into the remainder term. The modified remainder will be denoted by $\check{R}_T$.

\begin{lem}$\pt$ \label{mainlemnew}
For any $a \in D_H$, we have the following decomposition 
\begin{eqnarray} \label{maindeconew}
L_T(a) = L(a) +\frac{1}{T} H(a) + \frac{1}{T} K(a) + \frac{1}{T} \check{R}_T(a)
\end{eqnarray}
where $L$ and $H$ are given by \eqref{defLS} and \eqref{defHS} respectively,
\begin{equation} \label{defK2}
K(a)=-\frac{1}{2} \log \left( 1 + \frac{\vp(a) + \theta}{2\vp(a)} p_H \right),
\end{equation}
and the remainder term
\begin{align} \label{defW}
\check{R}_T(a) & =-\frac{1}{2}\log\left( 1 + \frac{(\vp(a) + \theta)(r_T(a)-p_H)}{(2+p_H)\vp(a) + \theta\delta_H} \right. \\
& + \left.
\frac{(\vp(a) + \theta)^2}{(\vp(a) - \theta) ((2+p_H)\vp(a) + \theta\delta_H)}e^{-2T\vp(a)}
\right). \nonumber
\end{align}
\end{lem}

\noindent{\bf Proof of Lemma \ref{mainlemnew}.}
It follows from \eqref{maindeco} that
\begin{align*}
L_T(a) &= L(a) +\frac{1}{T} H(a) + \frac{1}{T} K_T(a) + \frac{1}{T} R_T(a),\\
       &= L(a) + \frac{1}{T} H(a) + \frac{1}{T} K(a) + \frac{1}{T} \check{R}_T(a)
\end{align*}
where
$$
\check{R}_T(a)= K_T(a)- K(a) + R_T(a).
$$
Hence, denote $\vp = \sqrt{\theta^2 -2a}$, it is not hard to see that
\begin{align*}
\exp( -2 \check{R}_T(a))&=
\left( \frac{2\vp + (\vp + \theta)r_T(a)}{2\vp + (\vp + \theta)p_H} \right)
\left(1 + \frac{(\vp + \theta)^2}{(\vp - \theta)(2\vp+(\vp + \theta)r_T(a))}
e^{-2T\vp}\!\right),\\
&= \frac{2\vp + (\vp + \theta)r_T(a)}{2\vp + (\vp + \theta)p_H} +
\frac{(\vp + \theta)^2}{(\vp - \theta)(2\vp + (\vp + \theta)p_H)}e^{-2T\vp},\\
&= 1 + \frac{(\vp + \theta)(r_T(a)-p_H)}{2\vp + (\vp + \theta)p_H} +
\frac{(\vp + \theta)^2}{(\vp - \theta)(2\vp + (\vp + \theta)p_H)}e^{-2T\vp}.
\end{align*}
which ends the proof of Lemma \ref{mainlemnew}.\demend
\ \vspace{1ex}\\
Denote by $\Lambda_T$ the suitable approximation of the normalized cumulant generating function $L_T$ given by
\begin{equation} \label{def:laplace-approx-energie}
\Lambda_T(a) = L(a) + \frac{1}{T} H(a) + \frac{1}{T} K(a)
\end{equation}
The previous lemma enable us to write
\begin{eqnarray}
L_T(a) &=& \Lambda_T(a)+ \frac{1}{T} \check{R}_T(a).\label{def:local1}
\end{eqnarray}
One can observe that $\Lambda_T$ is an holomorph function in the domain $D_H$. In addition, there exists a unique $a_T$, which belongs to the interior of $D_H$ and converges to its border $a_H$, solution of the implicite equation
\begin{equation} \label{eqimpener}
\Lambda^{\prime}_T(a)=c.
\end{equation}
After some tedious but straightforward calculations, we can deduce from \eqref{eqimpener} that there exists a sequence $(a_k)$ such that, for any $p>0$ and $T$ large enough,
\begin{equation} \label{devaTener}
a_T =\sum_{k=0}^p \frac{a_k}{T^k}+{\cal O}\Bigl(\frac{1}{T^{p+1}}\Bigr)
\end{equation}
with
\begin{equation*}
a_0 = a_H, \qquad a_1 = -\frac{\theta \delta_H  }{  1 + 2\theta c \delta_H }, 
\end{equation*}
\begin{equation*}
a_2 =  \frac{\left( 2 \theta c \delta_H (4 + sin(\pi H) )  + 2 + sin(\pi H) \right)}{2(1 + 2 \theta c \delta_H)^3}.
\end{equation*}
Moreover, if $\varphi_T=\varphi(a_T)=\sqrt{\theta^2 - 2a_T}$, we also have for any $p>0$ and $T$ large enough,
\begin{equation} \label{devphiTener}
\varphi_T = \sum_{k=0}^p  \frac{\varphi_k}{T^k}+{\cal O}\Bigl(\frac{1}{T^{p+1}}\Bigr)
\end{equation}
with
\begin{equation*}
\varphi_0 = - \theta \delta_H, \qquad \varphi_1 = \frac{-1 }{1 + 2\theta c \delta_H},
\end{equation*}
\begin{equation*}
\varphi_2 = \frac{\left( 2 \theta c \delta_H  ( 5 + sin(\pi H) ) + 3 + sin(\pi H) \right)}{2\theta \delta_H(1 + 2 \theta c \delta_H)^3}.
\end{equation*}
Hereafter, we introduce the new probability measure
\begin{equation} \label{newprobaT}
\frac{d \p_T}{d \p }=\exp \Bigl( a_T S_T - T L_T(a_T) \Bigr)
\end{equation}
and we denote by $\E_T$ the expectation under $\p_T$. It clearly leads to the decomposition $\p(S_T \geq cT)=A_TB_T$
where
\begin{eqnarray}
\label{defAAener}
A_T &=& \exp \left( T L_T(a_T) - c T a_T  \right),\\
\label{defBBener}
B_T &=& \E_T \Bigl[\exp ( - a_T (S_T - c T)) \ind_{S_T \geq cT}\Bigr].
\end{eqnarray}
The proof now splits into two parts, the first one is devoted to the expansion of $A_T$ while the second one gives the expansion of $B_T$. It follows from \eqref{def:laplace-approx-energie}, \eqref{def:local1} and \eqref{defAAener} that
\begin{equation} \label{decoAAener}
A_T = \exp \Bigl( T (L(a_T) -ca_T)+H(a_T) + K(a_T) + \check{R}_T(a_T)\Bigr).
\end{equation}
We can deduce from the Taylor expansions of $a_T$ and $\varphi_T$ given by (\ref{devaTener}) and (\ref{devphiTener}) that
\begin{eqnarray*}
T(L(a_T) -c a_T)
&=& -\frac{T}{2}( \theta + \varphi_T +2ca_T),\\
&=& -T(c a_H -L(a_H))-\frac{\varphi_1}{2}-ca_1 -\frac{1}{2} \sum_{k=1}^p
\frac{\varphi_{k+1}+2ca_{k+1}}{T^k}+{\cal O}\Bigl(\frac{1}{T^{p+1}}\Bigr),\\
&=& -TI(c) +\frac{1}{2}-\frac{1}{2} \sum_{k=1}^p
\frac{\varphi_{k+1}+2ca_{k+1}}{T^k}+{\cal O}\Bigl(\frac{1}{T^{p+1}}\Bigr).
\end{eqnarray*}
Consequently, we obtain that for any $p>0$ and $T$ large enough,
\begin{equation} \label{devfinAAener1}
\exp \Bigl( T (L(a_T) -ca_T)\Bigr)=\exp(-TI(c)) \sqrt{e}\left[1 + \sum_{k=1}^{p}\frac{\alpha_{k}}{T^k}
+ {\cal O}\Bigl(\frac{1}{T^{p+1}}\Bigr)\right]
\end{equation}
where the coefficients $(\alpha_{k})$ may be explicitly calculated.
For example,
$$
\alpha_{1} =
\frac{-1}{4\theta \delta_H(1 + 2c\theta \delta_H)^2}
\left( 2\theta c \delta_H( 4 + sin(\pi H)) + 3 + sin(\pi H) \right).
$$
By the same way, we find that for any $p>0$ and $T$ large enough,
\begin{equation} \label{devfinAAener2}
\exp ( H(a_T) )=\sqrt{\frac{2\varphi_T}{\varphi_T-\theta}}=\sqrt{1 - \sin(\pi H)}
\left[1 + \sum_{k=1}^{p}\frac{\beta_{k}}{T^k}+ {\cal O}\Bigl(\frac{1}{T^{p+1}}\Bigr)\right]
\end{equation}
where the coefficients $(\beta_{k})$ may be explicitly calculated. For example,
$$
\beta_{1} = \frac{1 + \sin(\pi H)}{4 \theta \delta_H (1+2\theta c \delta_H)}.
$$
The expansions for $K(a_T)$ and $\check{R}_T(a_T)$ are much more tricky. On the one hand,
$$
\exp( K(a_T) ) = \sqrt{ \frac{ 2 \varphi_T }{ 2\varphi_T + (\varphi_T + \theta) p_H} }.
$$
One can observe that $2 \vp_0 + (\vp_0 + \theta)p_H=0$. Hence, multiplying the numerator and the denominator by $T$, we obtain that for any $p>0$ and $T$ large enough,
\begin{eqnarray} \label{devfinAAener3}
\exp( K(a_T) )
&=& \sqrt{ \frac{ 2 T \varphi_T }{ 2T \varphi_T + T(\varphi_T + \theta) p_H} }, \nonumber\\
&=& \sqrt{\theta T \delta_H(1-\delta_H)(1+2\theta c \delta_H)}\left[1 + \sum_{k=1}^{p}\frac{\gamma_{k}}{T^k} + {\cal O}\Bigl(\frac{1}{T^{p+1}}\Bigr)\right]
\end{eqnarray}
where, as before, the coefficients $(\gamma_{k})$ may be explicitly calculated. On the other hand, the remainder $\check{R}_T(a_T)= K_T(a_T)- K(a_T) + R_T(a_T)$. It is not hard to see that
\begin{align*}
\exp( -2 \check{R}_T(a_T))&= \frac{2\vp_T + (\vp_T +\theta)r_T(a_T)}{2\vp_T + (\vp_T + \theta)p_H} +
\frac{(\vp_T + \theta)^2}{(\vp_T-\theta)(2\vp_T +
(\vp_T +\theta)p_H)}\exp(-2T\vp_T),\\
&= \frac{2\vp_T + (\vp_T +\theta)r_T(a_T)}{2\vp_T + (\vp_T + \theta)p_H} +
{\cal O}(T\exp(2\theta T\delta_H)).
\end{align*}
Therefore,
$$
\exp( \check{R}_T(a_T))=\sqrt{\frac{2\vp_T + (\vp_T + \theta)p_H}{2\vp_T + (\vp_T +\theta)r_T(a_T)}}\Bigl(1+ {\cal O}(T\exp(2\theta T\delta_H))\Bigr).
$$
Recall that $r_T(a)=r_H(\vp(a)T/2)\exp(-T\vp(a))-1$. It is shown (equation \eqref{E-DEVRT} in Appendix A) that for any $p>0$ and $T$ large enough,
\begin{equation} \label{rTaTener}
r_T(a_T)=p_H+\frac{1}{\sin \pi H}\sum_{k=1}^{p}\frac{2^k r_{k}^H}{\vp_T^kT^k} + {\cal O}\Bigl(\frac{1}{T^{p+1}}\Bigr)
\end{equation}
where the coefficients $(r_{k}^H)$ may be explicitly calculated. For example,
$$
r_{1}^H=-\frac{(2H-1)^2}{4}.
$$
Consequently, we infer from (\ref{rTaTener}) that
$$
T(r_T(a_T)-p_H)=w_T(a_T)+ {\cal O}\Bigl(\frac{1}{T^{p}}\Bigr)
$$
where
$$
w_T(a_T)=\frac{1}{\sin (\pi H)}\sum_{k=1}^{p}\frac{2^k r_{k}^H}{\vp_T^kT^{k-1}}.
$$
If $\mu_T=T(2\vp_T + (\vp_T + \theta)p_H)$, we obtain that for any $p>0$ and $T$ large enough,
\begin{eqnarray}
\label{devfinAAener4}
\exp( \check{R}_T(a_T))&=&\sqrt{\frac{\mu_T}
{\mu_T +(\vp_T +\theta)T(r_T(a_T)-p_H)}}
\Bigl(1+ {\cal O}(T\exp(2\theta T\delta_H))\Bigr), \nonumber\\
&=&
\sqrt{\frac{1-(\sin( \pi H))^2}{1-(\sin( \pi H))^2+4r_1^H(\sin (\pi H))(1+2\theta c \delta_H)}}
\left[1 + \sum_{k=1}^{p}\frac{\delta_{k}}{T^k}
+ {\cal O}\Bigl(\frac{1}{T^{p+1}}\Bigr)\right]
\end{eqnarray}
where, as before, the coefficients $(\delta_{k})$ may be explicitly calculated. Putting together the four contributions (\ref{devfinAAener1}), (\ref{devfinAAener2}), (\ref{devfinAAener3}), and (\ref{devfinAAener4}), we find from (\ref{decoAAener}) that for any $p>0$ and $T$ large enough,
\begin{equation} \label{devfinAAener}
A_T\!=\!\exp(-TI(c)\!+\!R_H(c))\delta_H \sqrt{2 e \theta T \sin (\pi H)(1+2\theta c \delta_H)}
\!\left[1 \! + \!\sum_{k=1}^{p}\frac{\alpha_{k}}{T^k}
\!+\! {\cal O}\Bigl(\frac{1}{T^{p+1}}\Bigr)\right]
\end{equation}
where the coefficients $(\alpha_{k})$ may be explicitly calculated and
\begin{equation} \label{defRHinter}
R_H(c)=-\frac{1}{2} \log \left(1-\frac{(2H-1)^2\sin (\pi H)(1+2\theta c \delta_H)}{1-(\sin( \pi H))^2}\right).
\end{equation}
The rest of the proof concerns the expansion of $B_T$ which can be rewritten as
\begin{equation} \label{devBBener}
B_T = \E_T \Big[\exp(-a_T T U_T) \ind_{U_T\geq 0}\Big]
\end{equation}
where
$$
U_T = \frac{S_T -cT}{T}.
$$

\begin{lem} $\pt$ \label{lemBenerhard}
For all $c>- 1/(2\theta \delta_H)$, the distribution of $U_T$ under $\p_T$ converges, as $T$ goes to infinity, to the distribution of $\nu_H N^2-\gamma_H$ where $N$ stands for the standard $\cN(0,1)$ distribution,
\begin{equation} \label{defgammaH}
\gamma_H=c-L^{\prime}(a_H)=\frac{1+2\theta c \delta_H}{2\theta \delta_H},
\end{equation}
\begin{equation} \label{defnuH}
\nu_H=\frac{(1-(\sin (\pi H))^2)\gamma_H}{1-(\sin (\pi H))^2-(2H-1)^2\sin (\pi H)(1+2\theta c \delta_H)}.
\end{equation}
In other words, the limit of the characteristic function of $U_T$ under $\E_T$ is
\begin{equation} \label{limphienerhard}
\Phi(u) = \frac{\exp(-i \gamma_H u)}{\sqrt{1 - 2i \nu_H u }}.
\end{equation}
Moreover, there exists a sequence $(\beta_k)$ such that, for any $p>0$ and $T$ large enough,
\begin{equation} \label{devfinBBener}
B_T=\sum_{k=1}^p \frac{\beta_k}{T^k} + {\mathcal O}\Bigl( \frac{1}{T^{p+1}}\Bigr).
\end{equation}
The sequence $(\beta_k)$ only depends on the Taylor expansion of $a_T$ at the neighborhood of $a_H$ together with the derivatives of $L$ and $H$ evaluated at point $a_H$. They also depend on the Taylor expansion of $K_T$ and its derivatives at $a_H$. For example,
$$
\beta_1=\frac{1}{a_H\gamma_H \sqrt {2 \pi e}}\exp(Q_H(c)-R_H(c))
$$
where $R_H(c)$ is given by (\ref{defRHinter}) and
$$
Q_H(c)=\frac{(2H-1)^2\sin(\pi H)(1+2\theta c\delta_H)}{2(1-(\sin(\pi H))^2)}.
$$
\end{lem}

\noindent
{\bf Proof.} The proof of Lemma \ref{lemBenerhard} is given in Appendix B.\demend

\noindent{\bf Proof of Theorem \ref{theosldpener}, second part.} The expansions (\ref{devfinAAener}) and (\ref{devfinBBener}) imply (\ref{mainsldener3}), which completes the proof of Theorem \ref{theosldpener}. \demend

\subsection{Proof of Theorem \ref{theosldpenereg}.}

We shall now proceed to the proof of Theorem \ref{theosldpenereg} which essentially follows the same
lines as that of Theorem \ref{theosldpener}, second part. First of all, one can observe that if $c=- 1/(2\theta \delta_H)$, then we exactly have $a_c=a_H$. As in the proof of Theorem \ref{theosldpener}, there exists a unique $a_T$, which belongs to the interior of $D_H=]-\infty, a_H[$ and converges to its border $a_H$, solution of the implicite equation
\begin{equation} \label{eqimpenereg}
\Lambda^{\prime}_T(a)=c=-\frac{1}{2 \theta \delta_H}
\end{equation}
where $\Lambda_T$ is given by \eqref{def:laplace-approx-energie}. We deduce from \eqref{eqimpenereg} that
\begin{equation} \label{eqimpeneregbis}
T(\varphi_T + \theta \delta_H)^2= \frac{\theta (\varphi_T+ \theta p_H)}{c\varphi_T(\varphi_T-\theta)(2+p_H)}.
\end{equation}
Consequently, we infer from \eqref{eqimpenereg} and \eqref{eqimpeneregbis} that there exists a sequence $(a_k)$ such that, for any $p>0$ and $T$ large enough,
\begin{equation} \label{devaTphiTenereg}
a_T =\sum_{k=0}^{2p} \frac{a_k}{(\sqrt{T})^k}+{\cal O}\Bigl(\frac{1}{T^{p}\sqrt{T}}\Bigr)
\hspace{1cm}\text{and}\hspace{1cm}
\varphi_T = \sum_{k=0}^{2p}  \frac{\varphi_k}{(\sqrt{T})^k}+{\cal O}\Bigl(\frac{1}{T^{p}\sqrt{T}}\Bigr)
\end{equation}
with $a_0 = a_H$, $\varphi_0 = -\theta \delta_H$, $a_1= -(-\theta \delta_H)^{3/2}$, $\varphi_1=\sqrt{-\theta \delta_H}$,
\begin{align*}
a_2 &= - \frac{\theta \delta_H}{4}(1 + sin(\pi H)),\\
\varphi_2 &= - \frac{1}{4}(3 + sin(\pi H)).
\end{align*}
Furthermore, we have the decomposition $\p(S_T \geq cT)=A_TB_T$ where $A_T$ and $B_T$ are respectively given by \eqref{defAAener} and \eqref{defBBener}. Via the same lines as in the proof of the expansion \eqref{devfinAAener}, we find that for any $p>0$ and $T$ large enough,
\begin{equation} \label{devfinAAenereg}
A_T\!=\!\exp(-TI(c))(-\theta\delta_H e T)^{1/4}\sqrt{2 \delta_H \sin (\pi H)}\!\left[1 + \sum_{k=1}^{2p}\frac{\alpha_{k}}{(\sqrt{T})^k} + {\cal O}\Bigl(\frac{1}{T^{p}\sqrt{T}}\Bigr)\right]
\end{equation}
where the coefficients $(\alpha_{k})$ may be explicitly calculated. It still remains to give the expansion of $B_T$ which can be rewritten as
\begin{equation} \label{devBBenereg}
B_T = \E_T \Big[\exp(-a_T \sqrt{T} U_T) \ind_{U_T\geq 0}\Big]
\end{equation}
where
\begin{equation} \label{defUT}
U_T = \frac{S_T -cT}{\sqrt{T}}.
\end{equation}

\begin{lem}$\pt$ \label{lemBenereg}
For $c=- 1/(2\theta \delta_H)$, the distribution of $U_T$ under $\p_T$ converges, as $T$ goes to infinity, to the distribution of $\sigma_H N_1+ \nu_H (N_2^2-1)$ where $N_1$ and $N_2$ are two independent $\cN(0,1)$ random variables and
\begin{align}
\sigma_H^2 &=L^{\prime\prime }(a_H)= -\frac{1}{2(\theta \delta_H)^3},\\
\eta_H   &= \frac{1}{2(-\theta \delta_H)^{3/2}} \label{defetaH}.
\end{align}
In other words, the limit of the characteristic function of $U_T$ under $\E_T$ is
\begin{equation} \label{limphienereg}
\Phi(u) =  \frac{\exp\Bigl(-i \eta_H u - {\displaystyle \frac{u^2 \sigma_H^2}{2} } \Bigr)}{\sqrt{1 - 2i \eta_H u }}.
\end{equation}
Moreover, there exists a sequence $(\beta_k)$ such that, for any $p>0$ and $T$ large enough,
\begin{equation} \label{devfinBBenereg}
B_T=\sum_{k=1}^{2p} \frac{\beta_k}{(\sqrt{T})^k} + {\cal O}\Bigl(\frac{1}{T^{p}\sqrt{T}}\Bigr)
\end{equation}
where the sequence $(\beta_k)$ may be explicitly calculated. For example,
$$
\beta_1=\frac{1}{4 \pi a_H\eta_H}\exp\Bigl(-\frac{1}{4}\Bigr)\Gamma(\frac{1}{4}\Bigr).
$$
\end{lem}

\noindent{\bf Proof.} The proof of Lemma \ref{lemBenereg} is given in Appendix C.\demend

\noindent{\bf Proof of Theorem \ref{theosldpenereg}.} The expansions (\ref{devfinAAenereg}) and (\ref{devfinBBenereg}) imply (\ref{mainsldener4}), which completes the proof of Theorem \ref{theosldpenereg}.\demend

\section*{Appendix A: On the main asymptotic expansion.}
\renewcommand{\thesection}{\Alph{section}}
\renewcommand{\theequation}{\thesection.\arabic{equation}}
\setcounter{section}{1}
\setcounter{equation}{0}

We shall first prove the asymptotic expansion (\ref{maindeco}) of the normalized cumulant generating function $\mathcal{L}_T(a,b)$. This result was partially established by formula (5.12) of Kleptsyna and Le Breton \cite{KLB1}. By Girsanov's theorem, $\mathcal{L}_T(a,b)$ can be rewritten as \begin{eqnarray*}
\cL_T(a,b)&=&\frac{1}{T} \log \E\Bigl[\exp\Bigl(a\int_0^T \!Q_t\,dY_t+b S_T\Bigr)\Bigr],\\
&=&\frac{1}{T} \log \E_{\vp}\Bigl[\exp\Bigl((a+\theta-\vp)\int_0^TQ_t\,dY_t +\frac{1}{2}(2b-\theta^2+\vp^2)S_T\Bigr)\Bigr]
\end{eqnarray*}
for all $\vp \in \R$, where $\E_{\vp}$ stands for the expectation after the usual change of probability
\[
\frac{d\p_{\vp}}{d\p}=\exp\Bigl( (\vp-\theta)\int_0^TQ_t\,dY_t-\frac{1}{2} (\vp^2-\theta^2)S_T\Bigr).
\]
If $\theta^2-2b>0$, we can choose $\vp=\sqrt{\theta^2 -2b}$ and $\tau=\vp-(a+\theta)$ which leads to
\begin{eqnarray}
\label{eqL1}
\mathcal{L}_T(a,b)=\frac{1}{T}
\log \E_{\vp}\Bigl[\exp\Bigl(-\tau\!\int_0^TQ_t\,dY_t\Bigr)\Bigr].
\end{eqnarray}
By It\^o formula, we also have
\[
\int_0^TQ_t\,dY_t=\frac{1}{2}\Bigl(l_HY_T\int_0^Tt^{2H-1}\,dY_t-T\Bigr).
\]
Consequently, we obtain from (\ref{eqL1}) that
\begin{eqnarray}
\label{eqL2}
\mathcal{L}_T(a,b)=\frac{\tau}{2}+ \frac{1}{T}
\log \E_{\vp}\Bigl[\exp\Bigl(-\frac{\tau l_H}{2} Y_T\int_0^Tt^{2H-1}\,dY_t\Bigr)\Bigr].
\end{eqnarray}
Under the new probability $\p_{\vp}$, the pair $(Y_T,\int_0^Tt^{2H-1}\,dY_t)$ is Gaussian with mean zero and covariance matrix $\Gamma_T(\vp)$. Denote $I$ and $J$ the two matrices
\[
I=\left(\begin{array}{cc}
 1 & 0 \\
 0 & 1 \\
\end{array}\right)\hspace{1cm}
\mbox{and}\hspace{1cm}
J=\left(\begin{array}{cc}
 0 & 1 \\
 1 & 0 \\
\end{array}\right).
\]
As soon as the matrix
\[
M_T(a,b)=I+\frac{\tau l_H}{2}\Gamma_T^{1/2}(\vp)J\Gamma_T^{1/2}(\vp)
\]
is positive definite, we deduce from (\ref{eqL2}) together with standard calculus on the Gaussian distribution that
\begin{eqnarray} \label{eqL3}
\mathcal{L}_T(a,b)=\frac{\tau}{2}-\frac{1}{2T} \log \det(M_T(a,b)).
\end{eqnarray}
Furthermore, it was already proven by relation (5.12) of \cite{KLB1} that, if $\tau >0$
\begin{eqnarray} \label{eqM1}
\det(M_T(a,b))=\frac{1}{z_T}
\left[x_T\Bigl(1+\frac{\tau}{\vp}e^{\delta_T}
\sinh(\delta_T) \Bigr)^2 -
y_T \Bigl(1-\frac{\tau}{\vp} e^{\delta_T}
\cosh(\delta_T) \Bigr)^2 \right]
\end{eqnarray}
with $\delta_T=T\vp/2$, $x_T=I_{H-1}(\delta_T)I_{-H}(\delta_T)$, $y_T=I_{1-H}(\delta_T)I_{H}(\delta_T)$ and
\begin{eqnarray*}
z_T=x_T-y_T=\frac{4 \sin(\pi H)}{\pi \vp T}
\end{eqnarray*}
where $I_H$ is the modified Bessel function of the first kind. We refer the reader to \cite{Leb} Chapter 5 for the main properties of Bessel functions. Therefore, if $p_T=(x_T+y_T)/z_T$ and $r_T=2p_Te^{-T\vp}-1$, we deduce from (\ref{eqM1}) after some straightforward calculations that
\begin{eqnarray}
\det(M_T(a,b))&\!\!\!=\!\!\!&
\frac{(2\vp -\tau)^2}{4\vp^2}+p_T
\frac{\tau(2\vp-\tau)}{2\vp^2}e^{T\vp}+
\frac{\tau^2}{4\vp^2}e^{2T\vp}\nonumber \\
&\!\!\!=\!\!\!&\frac{\tau}{2\vp}e^{2T\vp}\left(
1\!+\!\frac{(2\vp-\tau)}{2\vp}r_T
\!+\!\frac{(2\vp -\tau)^2}{2\vp \tau}e^{-2T\vp}\right).
\label{eqM2}
\end{eqnarray}
Consequently, we infer from (\ref{eqL3}) and (\ref{eqM2}) that
\begin{eqnarray*}
\mathcal{L}_T(a,b)&=&-\frac{1}{2}(a+\theta +\vp)-\frac{1}{2T}
\log \left(\frac{\tau}{2\vp}\right) -\frac{1}{2T} \log \left(1+\frac{(2\vp-\tau)}{2\vp}r_T\right)\\
&-&\frac{1}{2T} \log \left(1+\frac{(2\vp -\tau)^2}{\tau(2\vp+r_T(2\vp-\tau))}e^{-2T\vp}\right).
\end{eqnarray*}
In order to complete the proof of Lemma \ref{mainlem}, it remains to show that the limiting domain $\Delta_H$ reduces to $\theta^2 -2b>0$ and $\sqrt{\theta^2 -2b}>\max(a+\theta;-\delta_H(a+\theta))$.
On the one hand, we already saw that our calculation is true as soon as $\theta^2 -2b>0$ and $\tau>0$ which can be rewritten as
$$
\vp>a+\theta.
$$
On the other hand, we also have the second constraint
\begin{equation} \label{eqC1}
1+\frac{(2\vp-\tau)}{2\vp}r_T>0
\end{equation}
leading to
\begin{equation} \label{eqC2}
\sqrt{\theta^2 -2b}>-\delta_H(a+\theta).
\end{equation}
As a matter of fact, it follows from the asymptotic expansion (5.11.10) of \cite{Leb} for the Bessel function $I_H$ that for all $z \in\C$ with $|z|$ large enough and $|\arg(z)|\leq \pi/2 -\delta$ where $\delta$ is an arbitrarily small positive number and for any $p>0$
\begin{eqnarray} \label{exprH}
r_H(z)=\frac{\exp(2z)}{sin(\pi H)} \left[1 + \sum_{k=1}^{p}\frac{r_{k}^H}{z^k} + {\cal O}\Bigl(\frac{1}{|z|^{p+1}}\Bigr)\right].
\end{eqnarray}
Moreover, the coefficients $(r_{k}^H)$ may be explicitly calculated. For example, one can check that $r_{1}^H=-(2H-1)^2/4$ and $r_{2}^H=(2H-1)^2(2H+1)(2H-3)/32$. In addition, all the coefficients $(r_{k}^H)$ vanish to zero if $H=1/2$. Consequently,
\begin{equation} \label{E-DEVRT}
r_T(a)=p_H+\frac{1}{\sin (\pi H)}\sum_{k=1}^{p}\frac{2^k r_{k}^H}{(\vp(a))^k T^k}
+ {\cal O}\Bigl(\frac{1}{T^{p+1}}\Bigr)
\end{equation}
with $p_H=(1- \sin(\pi H))/sin(\pi H)$. Hence, as $T$ tends to infinity, (\ref{eqC1}) reduces to $2\vp+(\vp+(a+\theta))p_H>0$ so $\vp(2+p_H)>-p_H(a+\theta)$. Finally, as $\delta_H=p_H/(2+p_H)$, it clearly implies (\ref{eqC2}) which completes the proof of Lemma \ref{mainlem}.\demend

\section*{Appendix B: On the characteristic functions.}

\renewcommand{\thesection}{\Alph{section}}
\renewcommand{\theequation}{\thesection.\arabic{equation}}
\setcounter{section}{2}
\setcounter{subsection}{0}
\setcounter{equation}{0}

\subsection{Proof of Lemma \ref{lemBenereasy}.}

If $\Phi_T$ denotes the characteristic function of $U_T$ under $\p_T$, it follows from \eqref{newprobstan} that
\begin{equation} \label{defcharenereasy}
\Phi_T(u) = \exp \left( -\frac{iuc\sqrt{T}}{\sigma_c} + T\left(L_T\Bigl(a_c+\frac{iu}{\sigma_c \sqrt{T}}\Bigr)-L_T(a_c) \right)\right).
\end{equation}
First of all, it is necessary to prove that for $T$ large enough, $\Phi_{T}$ belongs to $L^2(\R)$. One can observe that, in contrast with \cite{BR}, it is impossible here to make use of the Karhunen-Lo\`eve expansion of the process $(X_t)$.

\begin{lem}$\pt$ \label{lemcharenereasy}
For $T$ large enough, $\Phi_{T}$ belongs to $L^2(\R)$.
\end{lem}

\noindent{\bf Proof.} It is a direct consequence of Proposition \ref{energie-tf-l2} page \pageref{energie-tf-l2}.\demend

\noindent We shall now establish an asymptotic expansion for the characteristic function $\Phi_T$, similar to that of Lemma 7.1 of \cite{BR}.

\begin{lem}$\pt$ \label{expcharenereasy}
For any $p >0$, there exist integers $q(p)$, $r(p)$ and a sequence $(\varphi_{k,l}^H)$ independent of $p$, such that, for $T$ large enough
\begin{eqnarray} \label{devfincharenereasy}
\Phi_T(u)= \exp\Bigl(-\frac{u^2}{2}\Bigr)
\left[1+\frac{1}{\sqrt{T}}\sum_{k=0}^{2p}\sum_{l=k+1}^{q(p)}
\frac{\varphi_{k,l}^H u^l}{(\sqrt{T})^k}
+{\cal O}\Bigl(\frac{\max(1,|u|^{r(p)})}{T^{p+1}}\Bigr)\right]
\end{eqnarray}
and the remainder ${\cal O}$ is uniform as soon as $|u|\leq s T^{1/6}$ for some positive constant $s$.
\end{lem}

\noindent{\bf Proof.} It is rather easy to see that for all $k \in \N$, $R_T^{(k)}(a_c)={\cal O}(T^k\exp(-T/c))$. Hence, we infer from (\ref{maindeco}) together with (\ref{defLS}), (\ref{defHS}), (\ref{defKTS}) that for all $k \in \N$,
\begin{equation} \label{loglapenereasy}
L_T^{(k)}(a_c) = L^{(k)}(a_c) +\frac{1}{T}H^{(k)}(a_c) +\frac{1}{T}K_T^{(k)}(a_c)+{\cal O}(T^k\exp(-T/c)).
\end{equation}
Therefore, we find from (\ref{defcharenereasy}) and (\ref{loglapenereasy}) that for any $p >0$,
\begin{align*} \label{logfinenereasy}
\log \Phi_T(u) &= -\frac{u^2}{2} + T\sum_{k=3}^{2p+3}\left(\frac{iu}{\sigma_c\sqrt{T}}\right)^{k}\!  \frac{L^{(k)}(a_c)}{k!}\\
&+\sum_{k=1}^{2p+1} \left(\frac{iu}{\sigma_c \sqrt{T}}\right)^{k} \!  \frac{H^{(k)}(a_c)+K_T^{(k)}(a_c)}{k!}
+{\cal O}\Bigl(\frac{\max(1,u^{2p+4})}{T^{p+1}}\Bigr).
\end{align*}
We deduce the asymptotic expansion (\ref{devfincharenereasy}) by taking the exponential on both sides, remarking that, as soon as $|u|\leq s T^{1/6}$ some positive constant $s$, the quantity $u^l/(\sqrt{T})^k$ remains bounded in (\ref{devfincharenereasy}).\demend

\noindent{\bf Proof of Lemma \ref{lemBenereasy}.} It follows from Parseval's formula that $B_T$, given by (\ref{devBener}), can be rewritten as
\begin{equation} \label{ParBenereasy}
B_T=\frac{1}{2\pi a_c \sigma_c\sqrt{T}} \int_{\R}\left(1+\frac{iu}{a_c \sigma_c \sqrt{T}}\right)^{-1}\Phi_{T}(u) du.
\end{equation}
For some positive constant $s$, set $s_T=s T^{1/6}$. We can split $B_T=C_T + D_T$ with
\begin{eqnarray}
\label{defCTener1}
C_T&=& \frac{1}{2\pi a_c \sigma_c\sqrt{T}} \int_{|u| \leq s_T}\left(1+\frac{iu}{a_c \sigma_c \sqrt{T}}\right)^{-1}\Phi_{T}(u) du, \\
\label{defDTener2}
D_T&=& \frac{1}{2\pi a_c \sigma_c\sqrt{T}}  \int_{|u| > s_T}\left(1+\frac{iu}{ a_c \sigma_c \sqrt{T}}\right)^{-1}\Phi_{T}(u) du.
\end{eqnarray}
From now on, we claim that for some positive constant  $\nu$,
\begin{equation} \label{majDTener}
|D_T|={\cal O}(\exp(-\nu T^{1/3})).
\end{equation}
As a matter of fact, it follows from \eqref{defcharenereasy} that
\begin{equation*}
|\Phi_T(u) | \leq \exp \left(  T\left(L_T\Bigl(a_c+\frac{iu}{\sigma_c \sqrt{T}}\Bigr)-L_T(a_c) \right)\right).
\end{equation*}
We also deduce from \eqref{defLS} that $L(a_c) > 0$ and thus, using Proposition \ref{energie-prop-maj-r-a-d-t}, we find that
$$
|\Phi_{T}(u)| \leq \exp \left( - T L(a_c) \right) \, \exp
\left( -\frac{T|u|}{8 \sqrt{2}\vp(a_c)} 
\left( 1 + \frac{2|u|}{\vp^2(a_c)} \right)^{-3/4} \right)
$$
which leads to \eqref{majDTener}. Finally, we deduce (\ref{devfinBener}) from (\ref{devfincharenereasy}) and (\ref{defCTener1}) together with standard calculus on the $\cN(0,1)$ distribution.\demend

\subsection{Proof of Lemma \ref{lemBenerhard}.}

If $\Phi_T$ stands for the characteristic function of $U_T$ under $\p_T$, we have from (\ref{newprobaT})
\begin{equation} \label{defcharenerhard}
\Phi_T(u) = \exp \left( - iuc + T\left(L_T\Bigl(a_T+\frac{iu}{T}\Bigr)-L_T(a_T) \right)\right).
\end{equation}
As in the proof of Lemma \ref{lemBenereasy}, it follows from Proposition \ref{energie-tf-l2} page \pageref{energie-tf-l2} that for $T$ large enough, $\Phi_{T}$ belongs to $L^2(\R)$. We shall now propose an asymptotic expansion for $\Phi_T$, slightly different from that of Lemma 7.2 of \cite{BR}.
\begin{lem}$\pt$ \label{lemcharenerhard}
For any $p >0$, there exist integers $q(p)$, $r(p)$, $s(p)$ and a sequence $(\varphi_{k,l,m}^H)$ independent of $p$, such that, for $T$ large enough
\begin{equation*} \label{devfincharenerhard}
\Phi_T(u)= \Phi(u)\exp\left(-\frac{\sigma_H^2u^2}{2T}\right) \left[1+\!\sum_{k=1}^p\sum_{l=1}^{q(p)}\sum_{m=1}^{r(p)}\frac{\varphi^H_{k,l,m}u^l}{T^k(1-2i\nu_H u)^m}
+{\cal O}\Bigl(\frac{\max(1,|u|^{s(p)})}{T^{p+1}}\Bigr)\right]
\end{equation*}
where $\Phi$ is given by (\ref{limphienerhard}),
$$
\gamma_H=c-L^{\prime}(a_H)=\frac{1+2\theta c \delta_H}{2\theta \delta_H},\hspace{1cm} \sigma^2_H=L^{\prime\prime}(a_H)=-\frac{1}{2\theta^3\delta_H^3},
$$
and
$$
\nu_H=\frac{(1-(\sin \pi H)^2)\gamma_H}{1-(\sin \pi H)^2-(2H-1)^2\sin \pi H(1+2\theta c \delta_H)}.
$$
Moreover, the remainder ${\cal O}$ is uniform as soon as $|u|\!\leq\! s T^{2/3}$ for some positive constant $s$.
\end{lem}

\begin{rem}$\pt$
One can observe in this asymptotic expansion the limiting $\chi^2$ distribution $\Phi$ together with an independent centered Gaussian distribution with small variance $\sigma^2/T$.
\end{rem}

\noindent{\bf Proof.} First of all, we deduce from (\ref{maindeconew}) that
\begin{equation} \label{loglapenerhard}
L_T(a_T) = L(a_T) +\frac{1}{T}H(a_T) +\frac{1}{T}K(a_T)+\frac{1}{T}\check{R}_T(a_T).
\end{equation}
On the one hand, \eqref{defLS} implies that
\begin{equation*}
T\left(L\Bigl(a_T+\frac{iu}{T}\Bigr)-L(a_T) \right)=-\frac{T\varphi_T}{2}\left(\left(1-\frac{iub_T}{T}\right)^{1/2}-1\right)
\end{equation*}
with $b_T=2/\varphi_T^2$. Consequently, for all $p\geq 2$
\begin{equation*}
\exp \left(  T\left(L\left(a_T+\frac{iu}{T}\right)-L(a_T) \right)\right)=\exp\left(
\frac{iu\varphi_T b_T}{4}-\frac{T\varphi_T}{2} \sum_{k=2}^pl_k\Bigl(\frac{iub_T}{T}\Bigr)^k
+{\cal O}\Bigl(\frac{|u|^{p+1}}{T^{p}}\Bigr)\right)
\end{equation*}
where $l_k=-(2k)!/((2k-1)(2^k k!)^2)$ which leads to
\begin{equation}\label{devLhard}
\exp \left( -iuc+ T\left(L(a_T+\frac{iu}{T})-L(a_T) \right)\right)
\end{equation}
\vspace{-3ex}
\begin{equation*}
= \exp\Bigl(-iu\gamma_H - \frac{\sigma_H^2u^2}{2T} \Bigr)
\left[1+\sum_{k=1}^p\sum_{l=1}^{q(p)}\frac{\varphi^H_{k,l}u^l}{T^k}
+{\cal O}\Bigl(\frac{\max(1,|u|^{s(p)})}{T^{p+1}}\Bigr)\right].
\end{equation*}
On the other hand, we also have from (\ref{defHS}) that for all $p\geq 1$
\begin{align*}
\exp \left(  H\Bigl(a_T+\frac{iu}{T}\Bigr)-H(a_T) \right)&=\left(\frac{\varphi_T-\theta}
{\varphi_T-\theta(1-iub_T/T)^{-1/2}}\right)^{1/2},\\
&=\left(1- \Bigl(\frac{\theta}{\varphi_T-\theta}\Bigr)
\sum_{k=1}^ph_k
\Bigl(\frac{iub_T}{T}\Bigr)^k
+{\cal O}\Bigl(\frac{|u|^{p+1}}{T^{p+1}}\Bigr)\right)^{-1/2}
\end{align*}
with $h_k=(2k)!/(2^k k!)^2$. Hence,
\begin{equation}\label{devHhard}
\exp \left(  H\Bigl(a_T+\frac{iu}{T}\Bigr)-H(a_T) \right)
=\left[1+\sum_{k=1}^p\sum_{l=1}^{q(p)}
\frac{\psi^H_{k,l}u^l}{T^k}
+{\cal O}\Bigl(\frac{\max(1,|u|^{s(p)})}{T^{p+1}}\Bigr)\right].
\end{equation}
Furthermore, it follows from (\ref{defK2}) that for all $p\geq 1$
\begin{align*}
\exp \left(  K\Bigl(a_T+\frac{iu}{T}\Bigr)-K(a_T) \right)
&=\left(\frac{2\varphi_T + (\varphi_T+\theta)p_H}
{2\varphi_T + \varphi_Tp_H +\theta p_H
(1-iub_T/T)^{-1/2}}\right)^{1/2}, \\
&=\left(1+ \Bigl(\frac{\theta p_H T}{c_T}\Bigr)
\sum_{k=1}^ph_k
\Bigl(\frac{iub_T}{T}\Bigr)^k
+\frac{1}{c_T}{\cal O}\Bigl(\frac{|u|^{p+1}}{T^{p}}\Bigr)\right)^{-1/2}
\end{align*}
where $c_T=T(2\varphi_T+(\varphi_T+\theta)p_H)$. Therefore, if
$$
d_T(u)=1+\frac{iu\theta p_H b_T}{2c_T}
$$
we find that for all $p\geq 2$
\begin{equation}\label{deviniKhard}
\exp \left(  K\Bigl(a_T+\frac{iu}{T}\Bigr)-K(a_T) \right)
\end{equation}
\vspace{-3ex}
\begin{equation*}
=\frac{1}{\sqrt{d_T(u)}}\left(1+\Bigl(\frac{\theta p_H T}{c_Td_T(u)}\Bigr)
\sum_{k=2}^ph_k
\Bigl(\frac{iub_T}{T}\Bigr)^k
+\frac{1}{c_Td_T(u)}{\cal O}\Bigl(\frac{|u|^{p+1}}{T^{p}}\Bigr)
\right)^{-1/2}
\end{equation*}
One can easily check that as $T$ goes to infinity, the limits of $b_T$, $c_T$ and $d_T(u)$ are respectively given by $2/(\theta \delta_H)^2$, $-(2+p_H)/(1+2\theta c \delta_H)$ and $1-2i\gamma_H u$, where $\gamma_H$ is given by \eqref{defgammaH}. Then, we infer from \eqref{deviniKhard} that for all $p\geq 2$
\begin{equation}\label{devKhard}
\exp \left(  K\Bigl(a_T+\frac{iu}{T}\Bigr)-K(a_T) \right)
\end{equation}
\vspace{-3ex}
\begin{equation*}
=\frac{1}{\sqrt{1-2i\gamma_H u}}
\left[1+\sum_{k=1}^p\sum_{l=1}^{q(p)}
\sum_{m=1}^{r(p)}\frac{\Psi^H_{k,l,m}u^l}{T^k(1-2i\nu_H u)^m}
+{\cal O}\Bigl(\frac{\max(1,|u|^{s(p)})}{T^{p+1}}\Bigr)\right].
\end{equation*}
Now, in contrast with \cite{BR}, the remainder term $\check{R}_T$ plays a prominent role that can't be neglected. Let $\xi_T=T(\varphi_T+\theta)(r(a_T)-p_H)/c_T$ and
$$
\xi_T(u)=\frac{T}{c_T}\left(\varphi_T + \theta\Bigl(1-\frac{iub_T}{T}\Bigr)^{-1/2}\right)
\left(r_T\Bigl(a_T+\frac{iu}{T}\Bigr)-p_H\right).
$$
One can observe that $\xi_T$ and $\xi_T(u)$ share the same limit
$$
\lim_{T\rightarrow \infty}\xi_T(u)=\xi_H=\frac{2(1-\delta_H)(1+2\theta c \delta_H)r_1^H}{\delta_H (2+p_H) \sin(\pi H)}=-\frac{(2H-1)^2 \sin(\pi H) (1+2\theta c \delta_H)}{1-(\sin(\pi H))^2}.
$$
In addition, it follows from \eqref{defW} that
\begin{align*}
\exp \left( \check{R}_T(a_T) \right)
&=\left(\frac{c_T}
{c_T+c_T\xi_T}\right)^{1/2}\left[1+{\cal O}\Bigl(T\exp(2\theta T\delta_H)\Bigr)\right],\\
&=\left(1+\xi_T\right)^{-1/2}\left[1+{\cal O}\Bigl(T\exp(2\theta T\delta_H)\Bigr)\right].
\end{align*}
Moreover, we also have
\begin{equation*}
\exp \left(  \check{R}_T\Bigl(a_T+\frac{iu}{T}\Bigr) \right)=\left(\frac{d_T(u)+e_T(u)+\xi_T(u)}
{d_T(u)+e_T(u)}\right)^{-1/2}
\left[1+{\cal O}\Bigl(T\exp(2\theta T\delta_H)\Bigr)\right]
\end{equation*}
where
$$
e_T(u)=\frac{\theta p_H T}{c_T} \left(\Bigl(1-\frac{iub_T}{T}\Bigr)^{-1/2}-1-\frac{iub_T}{2T}\right).
$$
Therefore, via the same lines as in the proof of \eqref{devKhard}, we find that for all $p\geq 2$
\begin{equation}\label{devRhard}
\exp \left(  \check{R}_T\Bigl(a_T+\frac{iu}{T}\Bigr)-\check{R}_T(a_T) \right)
\end{equation}
\vspace{-3ex}
\begin{equation*}
=\frac{\sqrt{1-2i\gamma_H u}}{\sqrt{1-2i\nu_H u}}\left[1+\sum_{k=1}^p\sum_{l=1}^{q(p)}
\sum_{m=1}^{r(p)}\frac{\Phi^H_{k,l,m}u^l}{T^k(1-2i\nu_H u)^m}
+{\cal O}\Bigl(\frac{\max(1,|u|^{s(p)})}{T^{p+1}}\Bigr)\right]
\end{equation*}
with
\begin{equation*}
\nu_H=\frac{\gamma_H}{1+\xi_H}=\frac{(1-(\sin (\pi H))^2)\gamma_H}{1-(\sin (\pi H))^2-(2H-1)^2
\sin (\pi H)(1+2\theta c \delta_H)}.
\end{equation*}
Finally, Lemma \ref{lemcharenerhard} follows from the conjunction of \eqref{devLhard}, \eqref{devHhard}, \eqref{devKhard}, and \eqref{devRhard}.\demend

\noindent{\bf Proof of Lemma \ref{lemBenerhard}.} Via Parseval's formula, $B_T$, given by (\ref{devBBener}), can be rewritten as
\begin{equation} \label{ParBenerhard}
B_T=\frac{1}{2\pi T a_T}
\int_{\R}\left(1+\frac{iu}{T a_T}\right)^{-1}\Phi_{T}(u) du.
\end{equation}
Let $s_T>0$ such that $\sqrt{T}=o(s_T)$ as $T$ goes to infinity.
We can split $B_T=C_T + D_T$ where
\begin{eqnarray}
\label{defCTener1hard}
C_T&=&
\frac{1}{2\pi T a_T} \int_{|u| \leq s_T}
\left(1+\frac{iu}{T a_T}\right)^{-1}\Phi_{T}(u) du, \\
\label{defDTener2hard}
D_T&=&
\frac{1}{2\pi T a_T}  \int_{|u| > s_T}
\left(1+\frac{iu}{T  a_T}\right)^{-1}\Phi_{T}(u) du.
\end{eqnarray}
On the one hand, we find from Proposition \ref{energie-tf-l2},
and the fact that $x \mapsto x(1+x) ^{-3/4}$ is increasing that for some positive constant $\mu$, that
\begin{align*}
|D_T| ={ \cal O}\left(T(1+T^{3/2}) \exp\left(-\frac{ \mu s_T^2}{T} 
\Bigl(1 + \frac{s_T^2}{T^{2}}\Bigr)^{-3/4}\right) \right).
\end{align*}
It clearly leads to
$$| D_T |= {\cal O}(\exp(-\mu s_T^2/T)).$$
On the other hand, the asymptotic expansion for $C_T$, which immediately leads
to \eqref{devfinBBener}, follows from Lemma \ref{leminteg}, completing the proof
of Lemma \ref{lemBenerhard}.\demend

\subsection{Proof of Lemma \ref{lemBenereg}.}

The proof follows the same lines as the proof of Lemma \ref{lemBenerhard}. The most important difference is that the scale of Taylor expansion is in $\sqrt{T}$ instead of $T.$ Since $\Phi_T$ is the characteristic of $U_T$ defined by \eqref{defUT} under ${\mathbb P}_T$ defined by \eqref{newprobaT}, we have:
\begin{equation} \label{defcharenerlaure}
\Phi_T(u) = \exp \left( \frac{iu\sqrt{T}}{ 2 \theta \delta_H} +
T\left(L_T\Bigl(a_T+\frac{iu}{\sqrt{T}}\Bigr)-L_T(a_T) \right)\right).
\end{equation}
As in the proof of Lemma \ref{lemBenereasy}, it follows from Proposition \ref{energie-tf-l2} page \pageref{energie-tf-l2} that for $T$ large enough, $\Phi_{T}$ belongs to $L^2(\R)$. We shall now propose an asymptotic expansion for $\Phi_T$, slightly different from that of Lemma \ref{lemcharenerhard}.
\begin{lem}$\pt$ \label{lemcharenerlaure}
For any $p >0$, there exist integers $q(p)$, $r(p)$, $s(p)$ and a sequence $(\varphi_{k,l,m}^H)$ independent of $p$, such that, for $T$ large enough
\begin{equation*} \label{devfincharenerlaure}
\Phi_T(u)= \Phi(u) \left[1+\!\sum_{k=1}^p\sum_{l=1}^{q(p)}\sum_{m=1}^{r(p)}\frac{\varphi^H_{k,l,m}u^l}{\sqrt{T}^k(1-2i\nu_H u)^m}+{\cal O}\Bigl(\frac{\max(1,|u|^{s(p)})}{\sqrt{T}^{p+1}}\Bigr)\right]
\end{equation*}
where $\Phi$ is given by \eqref{limphienereg}. Moreover, the remainder ${\cal O}$ is uniform as soon as $|u|\!\leq\! s T^{1/6}$ for some positive constant $s$.
\end{lem}

\noindent{\bf Proof.} First of all, we deduce from (\ref{maindeconew}) that
\begin{equation}
\label{loglapenerlaure}
L_T(a_T) = L(a_T) +\frac{1}{T}H(a_T) +\frac{1}{T}K(a_T)+\frac{1}{T}\check{R}_T(a_T).
\end{equation}
On the one hand, \eqref{defLS} implies that
\begin{equation*}
T\left(L\Bigl(a_T+\frac{iu}{\sqrt{T}}\Bigr)-L(a_T) \right)=-\frac{T\varphi_T}{2}
\left(\left(1-\frac{iub_T}{\sqrt{T}}\right)^{1/2}-1\right)
\end{equation*}
with $b_T=2/\varphi_T^2$. Consequently, for all $p\geq 2$
\begin{equation*}
\exp \left(  T\left(L(a_T+\frac{iu}{\sqrt{T}})-L(a_T) \right)\right)=\exp\left(
\frac{iu\varphi_T b_T\sqrt{T}}{4}-\frac{T\varphi_T}{2} \sum_{k=2}^pl_k \Bigl(\frac{iub_T}{\sqrt{T}}\Bigr)^k +{\cal O}\Bigl(\frac{|u|^{p+1}}{\sqrt{T}^{p+1}}\Bigr)\right)
\end{equation*}
where $l_k=-(2k)!/((2k-1)(2^k k!)^2)$ which leads to
\begin{equation} \label{devLlaure}
\exp \left( \frac{iu\sqrt{T}}{2 \theta \delta_H}+ T\left(L(a_T+\frac{iu}{\sqrt{T}})-L(a_T) \right)\right)
\end{equation}
\vspace{-3ex}
\begin{equation*}
=\exp(-iu\eta_H - \frac{u^2\sigma_H^2}{2})
\left[1+\sum_{k=1}^p\sum_{l=1}^{q(p)}
\frac{\varphi^H_{k,l}u^l}{\sqrt{T}^k}
+{\cal O}\Bigl(\frac{\max(1,|u|^{s(p)})}{\sqrt{T}^{p+1}}\Bigr)\right].
\end{equation*}
On the other hand, we also have from (\ref{defHS}) that for all $p\geq 1$
\begin{align*}
\exp \left(  H\Bigl(a_T+\frac{iu}{\sqrt{T}}\Bigr)-H(a_T) \right)&=\left(\frac{\varphi_T-\theta}
{\varphi_T-\theta(1-iub_T/\sqrt{T})^{-1/2}}\right)^{1/2},\\
&=\left(1- \Bigl(\frac{\theta}{\varphi_T-\theta}\Bigr)
\sum_{k=1}^ph_k
\Bigl(\frac{iub_T}{\sqrt{T}}\Bigr)^k
+{\cal O}\Bigl(\frac{|u|^{p+1}}{\sqrt{T}^{p+1}}\Bigr)\right)^{-1/2}
\end{align*}
with $h_k=(2k)!/(2^k k!)^2$. Hence,
\begin{equation}\label{devHlaure}
\exp \left(  H\Bigl(a_T+\frac{iu}{\sqrt{T}}\Bigr)-H(a_T) \right)
=\left[1+\sum_{k=1}^p\sum_{l=1}^{q(p)}
\frac{\psi^H_{k,l}u^l}{\sqrt{T}^k}
+{\cal O}\Bigl(\frac{\max(1,|u|^{s(p)})}{\sqrt{T}^{p+1}}\Bigr)\right].
\end{equation}
Furthermore, it follows from (\ref{defK2}) that for all $p\geq 1$
\begin{align*}
\exp \left(  K\Bigl(a_T+\frac{iu}{\sqrt{T}}\Bigr)-K(a_T) \right)
&=\left(\frac{2\varphi_T + (\varphi_T+\theta)p_H}
{2\varphi_T + \varphi_Tp_H +\theta p_H
(1-iub_T/\sqrt{T})^{-1/2}}\right)^{1/2}, \\
&=\left(1+ \Bigl(\frac{\theta p_H \sqrt{T}}{c_T}\Bigr)
\sum_{k=1}^ph_k
\Bigl(\frac{iub_T}{\sqrt{T}}\Bigr)^k
+\frac{1}{c_T}{\cal O}\Bigl(\frac{|u|^{p+1}}{\sqrt{T}^{p}}\Bigr)\right)^{-1/2}
\end{align*}
where $c_T=\sqrt{T}(2\varphi_T+(\varphi_T+\theta)p_H)$. Therefore, if
$$
d_T(u)=1+\frac{iu\theta p_H b_T}{2c_T},
$$
we find that for all $p\geq 2$
\begin{equation} \label{deviniKlaure}
\exp \left(  K\Bigl(a_T+\frac{iu}{\sqrt{T}}\Bigr)-K(a_T) \right)
\end{equation}
\vspace{-3ex}
\begin{equation*}
=\frac{1}{\sqrt{d_T(u)}}\left(1+\Bigl(\frac{\theta p_H T}{c_Td_T(u)}\Bigr)
\sum_{k=2}^ph_k
\Bigl(\frac{iub_T}{\sqrt{T}}\Bigr)^k
+\frac{1}{c_Td_T(u)}{\cal O}\Bigl(\frac{|u|^{p+1}}{\sqrt{T}^{p}}\Bigr)
\right)^{-1/2}
\end{equation*}
One can easily check that as $T$ goes to infinity, the limits of $b_T$, $c_T$ and $d_T(u)$ are respectively given by $2/(\theta \delta_H)^2$, $(2+p_H)\sqrt{ - \theta \delta_H}$ and $1-2i\eta_H u$, where $\eta_H$ is given by \eqref{defetaH}. Then, we infer from \eqref{deviniKlaure} that for all $p\geq 2$
\begin{equation} \label{devKlaure}
\exp \left(  K\Bigl(a_T+\frac{iu}{\sqrt{T}}\Bigr)-K(a_T) \right)
\end{equation}
\vspace{-3ex}
\begin{equation*}
=\frac{1}{\sqrt{1-2i\eta_H u}}
\left[1+\sum_{k=1}^p\sum_{l=1}^{q(p)}
\sum_{m=1}^{r(p)}\frac{\Psi^H_{k,l,m}u^l}{\sqrt{T}^k(1-2i\nu_H u)^m}
+{\cal O}\Bigl(\frac{\max(1,|u|^{s(p)})}{\sqrt{T}^{p+1}}\Bigr)\right].
\end{equation*}
Now, in contrast with \cite{BR}, the remainder term $\check{R}_T$ plays a prominent role that can't be neglected. Let $\xi_T=\sqrt{T}(\varphi_T+\theta)(r(a_T)-p_H)/c_T$ and
$$
\xi_T(u)=\frac{\sqrt{T}}{c_T}\left(\varphi_T +\theta\Bigl(1-\frac{iub_T}{\sqrt{T}}\Bigr)^{-1/2}\right)
\left(r_T\Bigl(a_T+\frac{iu}{\sqrt{T}}\Bigr)-p_H\right).
$$
One can observe that $\xi_T$ and $\xi_T(u)$ share the same limit
$$
\lim_{T\rightarrow \infty}\xi_T(u)=0.
$$
In addition, it follows from \eqref{defW} that
\begin{align*}
\exp \left( \check{R}_T(a_T) \right)
&=\left(\frac{c_T}
{c_T+c_T\xi_T}\right)^{1/2}\left[1+{\cal O}\Bigl(\sqrt{T}\exp(2\theta T\delta_H)\Bigr)\right],\\
&=\left(1+\xi_T\right)^{-1/2}\left[1+{\cal O}\Bigl(\sqrt{T}\exp(2\theta T\delta_H)\Bigr)\right].
\end{align*}
Moreover, we also have
\begin{equation*}
\exp \left(  \check{R}_T\Bigl(a_T+\frac{iu}{T}\Bigr) \right)=\left(\frac{d_T(u)+e_T(u)+\xi_T(u)}
{d_T(u)+e_T(u)}\right)^{-1/2}
\left[1+{\cal O}\Bigl(\sqrt{T}\exp(2\theta T\delta_H)\Bigr)\right]
\end{equation*}
where
$$
e_T(u)=\frac{\theta p_H T}{c_T} \left(\Bigl(1-\frac{iub_T}{\sqrt{T}}\Bigr)^{-1/2}-1
-\frac{iub_T}{2\sqrt{T}}\right).
$$
Therefore, via the same lines as in the proof of \eqref{devKlaure}, we find that for all $p\geq 2$
\begin{equation} \label{devRlaure}
\exp \left(  \check{R}_T\Bigl(a_T+\frac{iu}{T}\Bigr)-\check{R}_T(a_T) \right)
\end{equation}
\vspace{-3ex}
\begin{equation*}
=1+\sum_{k=1}^p\sum_{l=1}^{q(p)}
\sum_{m=1}^{r(p)}\frac{\Phi^H_{k,l,m}u^l}{\sqrt{T}^k(1-2i\nu_H u)^m}
+{\cal O}\Bigl(\frac{\max(1,|u|^{s(p)})}{\sqrt{T}^{p+1}}\Bigr).
\end{equation*}

Finally, Lemma \ref{lemcharenerlaure} follows from the conjunction of \eqref{devLlaure}, \eqref{devHlaure}, \eqref{devKlaure}, and \eqref{devRlaure}.\demend

\noindent{\bf Proof of Lemma \ref{lemBenereg}.} Via Parseval's formula, $B_T$, given by (\ref{devBBener}), can be rewritten as
\begin{equation} \label{ParBenereg}
B_T=\frac{1}{2\pi T a_T} \int_{\R}\left(1+\frac{iu}{T a_T}\right)^{-1}\Phi_{T}(u) du.
\end{equation}
For some positive constant $s$, set $s_T=s \, T^{1/6}$. We can split $B_T=C_T + D_T$ where
\begin{eqnarray}
\label{defCTener1eg}
C_T&=&
\frac{1}{2\pi T a_T} \int_{|u| \leq s_T}
\left(1+\frac{iu}{T a_T}\right)^{-1}\Phi_{T}(u) du, \\
\label{defDTener2eg}
D_T&=&
\frac{1}{2\pi T a_T}  \int_{|u| > s_T}
\left(1+\frac{iu}{T  a_T}\right)^{-1}\Phi_{T}(u) du.
\end{eqnarray}
On the one hand, we find from Proposition \ref{energie-tf-l2} that for some positive constant $\mu$,
$$
| D_T |= {\cal O}(\exp(-\mu s_T^2/T)).
$$
On the other hand, the asymptotic expansion for $C_T$, which immediately leads to \eqref{devfinBBener}, following the same arguments as those of Bercu and Rouault \cite[page 18]{BR}.\demend


\section*{Appendix C: Technical results.}
\renewcommand{\thesection}{\Alph{section}}
\renewcommand{\theequation}{\thesection.\arabic{equation}}
\setcounter{section}{3}
\setcounter{subsection}{0}
\setcounter{equation}{0}


\subsection{Statement of the results}

The main interest of the decomposition \eqref{maindeconew} is given by the two following results. They show us that the different functions we deal with are holomorphs and that the behaviour of the remainder is negligeable in our calculations.

\begin{prop}\label{energie-prop-holo}
Denote
$$
\cD_{\Delta} = \{ a\in {\mathbb C},~~ Re(a) < a_H \}
$$
and, for $\varepsilon>0,$
$$
\cD_1 = \left\{ a \in {\mathbb C},~  Re(a) < a_H - \varepsilon(2+\varepsilon)\frac{\theta^2 \delta_H^2}{2} \right\}.
$$
Then, for $T$ large enough, we have the following assertions.
\begin{enumerate}
\item[a)] The functions $\vp,$  $L_T,$ $L$, $H$, $K$ and $\check{R}_T$ have analytic extensions to $\cD_{\Delta}$.
\item[b)] The function $(a,T) \mapsto \check{R}_T(a)$ is $C^{\infty}$ on $\cD_{\Delta} \times [T_{\Delta},+ \infty[,$ for $T_{\Delta}$ depending only on $H$ and $\theta$.
\item[c)] For $\varepsilon>0$, and for all $a \in {\mathcal D}_1$, 
\begin{align*}
\sqrt{\frac{\delta_H}{ (2 +p_H)(\delta_H +1)}} \leq \Bigl| \exp(H(a)+K(a))\Bigr|   \leq \frac{4}{\sqrt{2+p_H}} \sqrt{\frac{1+ \varepsilon}{\varepsilon}}.
\end{align*}
\item[d)] There exists a constant $C$ depending only on $\theta$ and $H$ such that for $T$ large enough and for all $a \in \cD_1$,
\begin{align*}
\sqrt{ \frac{1}{2} -\frac{C}{ T^2 \varepsilon (\delta_H \theta)^4}}\leq \left|\exp(-\check{R}_T (a))\right|  \leq C \sqrt{ 1 + \frac{1}{T} + \frac{1}{T \varepsilon}}.
\end{align*}
\end{enumerate}
\end{prop}

\begin{prop}\label{energie-tf-l2}
For $T$ large enough, for $\varepsilon > 2C/( T^2 (\delta_H \theta)^4),$ $a \in \cD_1\cap {\mathbb R}$ and $u \in {\mathbb R},$
\begin{align*}
\Bigl| \exp( T ( L_T(a+iu) -L_T(a))) \Bigr| \leq C (1+T^{3/2}) \exp
\left( -\frac{T|u|}{8 \sqrt{2}\varphi(a)}  \left( 1 + \frac{2|u|}{\varphi^2(a)} \right)^{-3/4}\right)
\end{align*}
and the map
\begin{align*}
u \mapsto \exp( T \left( L_T( a + iu) - L_T(a) \right) )
\end{align*}
belongs to $L^2(\R).$
\end{prop}

\begin{prop}\label{energie-prop-maj-r-a-d-t}
As $T$ goes to infinity and $a \in {\mathbb R}$ such that $a < a_H$, we have
\begin{align*}
\exp( -\check{R}_T(a) ) = {\cal O} \left( \max \left(1 ; \frac{-1}{T(\varphi(a)+\delta_H \theta)}\right) \right).
\end{align*}
\end{prop}

\subsection{Proofs of the results}

We shall denote the principal determination of the logarithm defined on $\C \setminus ]-\infty, 0]$ by
\begin{align*}
\loga{z} = \log |z| + i Arg(z),
\end{align*}
where
\begin{eqnarray*}
Arg(z)
=& \arcsin \left[ \frac{ Im(z)}{|z|} \right]       &\mbox{~~if~~} Re(z) \geq 0,\\
=& \arccos \left[ \frac{Re(z)}{|z|} \right]        &\mbox{~~if~~} Re(z) <0,~~Im(z)>0,\\
=&  -\arccos \left[ \frac{Re(z)}{|z|} \right]   &\mbox{~~if~~} Re(z) <0,~~Im(z)<0,
\end{eqnarray*}

\noindent{\bf Proof of Proposition~\ref{energie-prop-holo}:} Since $T \mapsto S_T$ is a positive increasing process,
then for  $a <a_H$
\begin{align*}
\esp{}{ \exp (a S_T) } \leq \lim_{ T \rightarrow + \infty} \exp( T L_T(a)) < \infty.
\end{align*}
Lebesgue dominated theorem 
 yields that $a \mapsto \esp{}{\exp(aS_T)}$ has an analytic extension to $\{a \in {\mathbb C},~~ Re(a) < a_H.\}.$
In order  to prove Proposition  \ref{energie-prop-holo}, we have to obtain 
the same result for $\vp$, $L,$ $H,$ $K$ and $\check{R}_T.$ The proof is split in steps. First, we study the function $\vp.$

\begin{lem}\label{lem-energie-fi}
The function $\vp$ has an analytic extension on $\{a\in {\mathbb C},~~ Re(a) < a_H\},$ still denoted by $\vp$ such that $Arg(\varphi) \in ]-\frac{\pi}{4},\frac{\pi}{4},[,$ $Re(\varphi) \in ]-\delta_H \theta, + \infty[,$ $Im(\varphi)(a)$ vanishes if and only if $Im(a)=0,$ and for $\varepsilon > 0,$
\begin{align*}
\inf_{ a \in \cD_1 } \left\{ Re(\vp(a)) \right\} \, > \, - \theta \delta_H(1+\varepsilon).
\end{align*}
For all $a >a_H $ and $ z \in {\mathbb R},$
\begin{align}
Re( \vp(a+iz) - \vp(a) )
&\geq  \frac{|z|}{2 \sqrt{2}\varphi(a)} \left( 1 + \frac{2|z|}{\varphi^2(a)} \right)^{-3/4}. \label{E-MAJREF}
\end{align}
\end{lem}

\noindent{\bf Proof of Lemma \ref{lem-energie-fi}} The properties of $\vp$ relies on the properties of the analytic function defined on $Re(z)>0$ by
\begin{align*}
\sqrt{1+ z}= \sqrt{|1+z|} \expo{ \frac{i}{2} \, Arg(1+z) },
\end{align*}
and the fact that $\varphi(a)= -\delta_H \theta \sqrt{1+ \frac{2}{\delta_H^2 \theta^2} (a_H-a)}.$\\
Since  for $Re(z)>0,$ $Arg(z+1)$ belongs to $]-\frac{\pi}{2},\frac{\pi}{2},[,$ then $Arg \sqrt{1+ z}$ belongs to $]-\frac{\pi}{4},\frac{\pi}{4},[.$ Its imaginary part is
\begin{align*}
Im(\sqrt{1+z})= \frac{ \sqrt{ |1+z| - Re(1+z)}}{\sqrt{2}} sign(Im(z)).
\end{align*}
Its real part is given by
\begin{align*}
Re(\sqrt{1+z})= \frac{ \sqrt{ |1+z| +Re(1+z)}}{\sqrt{2}},
\end{align*}
is an increasing function of $Im(z)$ then fulfills
\begin{align*}
Re(\sqrt{1+z}) -1 \geq \frac{ Re(z)}{1 + \sqrt{Re(1+z)}},
\end{align*}
and for all $z,$ such that $Re(z) >\varepsilon(2 + \varepsilon),$
$Re(\sqrt{1+z})-1 \geq \varepsilon$. 
Inequality \eqref{E-MAJREF} is a consequence of the fact that for $x \geq 0$,
$$
\sqrt{1 + \sqrt{1 + x}} - \sqrt{2} \geq \frac{x}{4 \sqrt{2} ( 1 + x )^{\frac{3}{4}}}.
$$
\demend
Hereafter, we study $H$ and $K$. Observe that for $a \in ]-\infty,0],$ $H(a)= \tilde{H}(\vp(a))$ and
$K(a)= \tilde{K}(\vp(a))$ where for $z\in {\mathbb C},$ $Re(z)>-\delta_H \theta$,
\begin{align*}
\tilde{H}(z) &=-\frac{1}{2} \loga{ \frac{z-\theta}{2z} },\\
\tilde{K}(z) &=-\frac{1}{2} \loga{ 1 + \frac{z + \theta}{2z} p_H }.
\end{align*}
Then, the expected properties on $H$ and $K$ are some consequences of Lemma \ref{lem-energie-fi} and the same properties of $\tilde{H}$ and $\tilde{K}$ on $\{z\in{\mathbb C}	,~~ Re(z) > -\delta_H \theta\}$ instead of $\{a\in{\mathbb C},~~ R e(a) <a_H\}.$

\begin{lem} \label{energie-lem-tilde-H}
The functions $\tilde{H}$ and $\tilde{K}$ admit analytical extensions on $\{z,~~Re(z) > - \delta_H \theta\}$. Moreover, for all $\varepsilon>0,$ we denote by
$$
\cD_2 = \left\{ z\in{\mathbb C},~~ Re(z) >-\delta_H \theta (1+ \varepsilon),~~|Arg(z)|\leq  \frac{\pi}{4} \right\},
$$
and we have for $z \in {\mathcal D}_2$
\begin{align*}
\sqrt{\frac{\delta_H}{ (2 +p_H)(\delta_H +1)}}
\leq   \left| \exp\Bigl( \tilde{H}(z)+\tilde{K}(z) \Bigr) \right|
\leq  \frac{4}{\sqrt{2+p_H}} \sqrt{\frac{1+\varepsilon}{\varepsilon}}.
\end{align*}
\end{lem}

\noindent{\bf Proof of lemma \ref{energie-lem-tilde-H}~:} Using the expression of the functions $\tilde{H}$ and $\tilde{K}$, it is easy to prove the analytical extensions of $\tilde{H}$ and $\tilde{K}$ on $\{z\in {\mathbb C},~Re(z) > - \delta_H \theta\}$.  For $\tilde{K} $ it is a consequence of the fact that for $z\in {\mathbb C},$ $Re(z) >-\delta_H \theta,$ 
\begin{align*}
Re( 1 + \frac{z+ \theta}{2z} p_H) \geq  1 + \frac{p_H}{2} + \frac{\theta p_H}{ 2Re(z)} > 1 + \frac{p_H}{2} + \frac{ p_H}{ -2\delta_H} =0.
\end{align*}
It remains to prove the inequality states in the lemma but, using the fact that,
\begin{align*}
\tilde{H}(z) + \tilde{K}(z)
&= -\frac{1}{2}\left( \loga{z-\theta} + \loga{2z + (z + \theta)p_H }   - 2\loga{2z} \right),\\
&= -\frac{1}{2}\left( \loga{z-\theta} + \loga{(2+p_H) z + \theta p_H } - 2\loga{2z} \right),\\
&= -\frac{1}{2}\left( \loga{z-\theta} + \loga{2+p_H}                    + \loga{z + \frac{\theta p_H}{2+p_H} } - 2\loga{2z} \right),\\
&= -\frac{1}{2}\left( \loga{2+p_H} + \loga{z-\theta} + \loga{z+ \delta_H \theta}- \loga{4} -2\loga{z} \right).
\end{align*}
thus
$$
\exp\Bigl( \tilde{H}(z)+\tilde{K}(z) \Bigr) = \left(  \frac{(2 + p_H)(z-\theta)(z+\theta \delta_H)}{4z^2}  \right)^{-1/2},
$$
where $\sqrt{z}= \sqrt{|z|} \exp(i Arg(z)/2)$.
Since for $z \in {\mathbb C}$,  $|Arg(z)| \leq \frac{\pi}{4},$ $|Im(z) | \leq Re(z),$ if moreover we have $Re(z) > -\delta_H \theta (1+\varepsilon),$
then
\begin{align*}
1 \geq \left|\frac{z+ \delta_H \theta}{z}\right| \geq \frac{Re(z) + \delta_H \theta}{ Re(z)} \geq \frac{\varepsilon}{ ( 1 +\varepsilon)}
\end{align*}
and
\begin{align*}
\frac{1+ \delta_H}{\delta_H}\geq \left|\frac{z-\theta}{ z}\right|\geq \frac{ Re(z) - \theta}{ 2 Re(z)} \geq 1.
\end{align*}
We have used the fact that the function $x\mapsto \frac{x+\delta_H \theta}{x}$ is increasing and $x \mapsto \frac{x-\theta}{x}$ is decreasing. The third point of proposition \ref{energie-prop-holo} is then proved.
\demend
Now, we focus on $\check{R}_T.$ First, observe that for $a \in {\mathbb C},$ such that  $Re(a) <0,$ $\check{R}_T(a)= \tilde{\check{R}}_T(\varphi(a)),$ where for $z \in {\mathbb C},$ such that   $Re(z) >-\delta_H \theta,$
\begin{align} \label{-energie-def-tilde-r}
\tilde{\check{R}}_T(z)
&= -\frac{1}{2} \loga{ 1 + \frac{(z+\theta)(\tilde{r}_T(z)-p_H)}{(2+p_H)(z+ \delta_H \theta)} + \frac{(z+ \theta)^2}{ (2+p_H)(z-\theta)(z+\delta_H \theta)}e^{-2Tz} } \notag\\
&= -\frac{1}{2} \loga{ z(2+p_H)(z+\delta_H \theta) + z(z+\theta)(\tilde{r}_T(z)-p_H) +
\frac{z(z+ \theta)^2}{(z-\theta)}e^{-2Tz} }\\
&\quad+ \frac{1}{2} \loga{ z(2+p_H)(z+\delta_H \theta) }.\notag\\
&= -\frac{1}{2} \left( \loga{ \tilde{\check{R}}_{1,T}(z) } - \loga{(2+p_H)z} - \loga{z+\delta_H \theta} \right).\label{E-DECRCT}
\end{align}
where
\begin{align*}
\tilde{r}_T(z)             &= r_H \left( \frac{Tz}{2} \right) e^{-Tz} -1.\\
\tilde{\check{R}}_{1,T}(z) &= z(2+p_H)(z+\delta_H \theta) + z(z+\theta)(\tilde{r}_T(z)-p_H) +
\frac{z(z+ \theta)^2}{(z-\theta)}e^{-2Tz}.
\end{align*}
The properties of $\check{R}_T$ are some consequences of lemma \ref{lem-energie-fi} and the following lemma.

\begin{lem}\label{energie-lem-ppte-r}
Denote
$$
\cD_3 = \left\{ z,~ Re(z) > - \delta_H \theta,~~ Arg(z) \in ]-\frac{\pi}{4},\frac{\pi}{4}[ \right\}.
$$
For $T$ large enough and for all $z \in \cD_3$, $\tilde{\check{R}}_{1,T}(z) \in {\mathbb C} \setminus ]- \infty,0].$
\end{lem}

In fact, lemma \ref{energie-lem-ppte-r} and decomposition \eqref{E-DECRCT} give us an analytical extension of $\tilde{\check{R}}_{T}$ on $\cD_4$. Moreover, there exists $T_4$ depending only on $H$ and $\theta$  such that the function
\begin{align*}
(z ,T) \mapsto \frac{\tilde{\check{R}}_{1,T}(z)}{(2+p_H) z(z +\delta_H \theta)}
\end{align*}
is $C^{\infty}$ and never vanishes on $\{ (z,T)~,~z \in \cD_4, \, T > T_4 \}$ and is $C^{\infty}$ with respect to $(z,T).$ Since, $\check{R}_T = \tilde{\check{R}}_{T}(\vp),$ then $(a,T) \mapsto \check{R}_T(a) $ is $C^{\infty}$ on $T > T_4$ and $\{a,~~Re(a)< a_H\}.$ The second point of proposition \ref{energie-prop-holo} is then proved.
\demend

\noindent{\bf Proof of Lemma \ref{energie-lem-ppte-r}.}
We recall that for $z\in {\mathbb C}$  such that $Arg(z) \in ]-\frac{\pi}{4},\frac{\pi}{4},[$
\begin{align*}
r_H(z)-1=\frac{\exp(2z)}{\sin(\pi H)} \left( 1 - \frac{(2H-1)^2}{4z} + \frac{1}{z^2} \sin(\pi H) F(z) \right)
\end{align*}
where $F$ is a continuous bounded function. Then, for  $z \in \cD_3$,
\begin{align*}
\tilde{r}_T(z) -p_H = \frac{r_1}{ Tz} + \frac{1}{T^2z^2} \, F(Tz).
\end{align*}
We can exhibit a polynomial term from $\tilde{\check{R}}_{1,T}(z)$
\begin{align}
\tilde{\check{R}}_{1,T}(z)
&= z(2+p_H)(z+\delta_H \theta) + z(z+\theta) \left( \frac{r_1}{ Tz} + \frac{1}{T^2z^2} \, F(Tz) \right) + \frac{z(z+ \theta)^2}{(z-\theta)}e^{-2Tz},\notag\\
&= z(2+p_H)(z+\delta_H \theta) + (z+\theta)\frac{r_1}{ T}  + \left(\frac{z+\theta}{T^2z} \, F(Tz)  + \frac{z(z+ \theta)^2}{(z-\theta)}e^{-2Tz}\right),\notag\\
&=z(2+p_H)(z+\delta_H \theta) + (z+\theta)\frac{r_1}{ T} + \frac{z+\theta}{T^2z} \, \left( F(Tz)  + \frac{T^2z^2(z+ \theta)}{(z-\theta)}e^{-2Tz}\right). \label{E-DEFR1CT}
\end{align}
Let us denote
\begin{align}
\tilde{P}_T(z)
&= z(2+p_H)(z+\delta_H \theta) + (z+\theta)\frac{r_1}{ T}, \label{E-DEFP}\\
C
&= \frac{1- \delta_H}{\delta_H}\sup_{ \{ z \in \cD_3 \}} \left\{ \left| F(Tz)+\frac{(z+ \theta)}{ (z-\theta)} T^2 z^2e^{-2zT}\right| \right\} < +\infty \notag.
\end{align}
Observe that for $z \in {\mathcal D}_3,$ on the one hand, 
\begin{align*}
\left| Im(\tilde{P}_T(z)) \right|&=\left|Im(z)\right| \left[ (2+p_H) (2 Re(z) + \delta_H \theta) +\frac{r_1}{T} \right]\\
&\geq \left|Im(z)\right|\left[ \frac{(2+p_H)}{2}(- \delta_H \theta) + \frac{r_1}{T}\right] .
\end{align*}
Then, for $z \in {\mathbb C},$ such that 
\begin{align*}
Im(z) > \frac{4}{ 2 +p_H} \frac{1}{- \delta_H \theta} \left[ \frac{C}{T^2} + \frac{-r_1}{T}\right],
\end{align*}
 $|Im(\tilde{\check{R}}_{1,T}(z))| >0.$ 
On the other hand, 
\begin{align*}
Re (\tilde{P}_T(z) )&= (2+p_H) Re(z)(Re(z) + \delta_H \theta) + \frac{\theta r_1}{T} Re(z) - (2+p_H) Im(z)^2,\\
&\geq -\delta_H \theta^2 \frac{r_1}{T} -(2+p_H) Im(z)^2,
\end{align*}
Then, for $T$ large enough, for all $z \in {\mathcal D}_3$ such that 
if \begin{align*}
Im(z) > \frac{4}{ 2 +p_H} \frac{1}{- \delta_H \theta} \left[ \frac{C}{T^2} + \frac{-r_1}{T}\right],
\end{align*}
then $Re( \tilde{\check{R}}_{1,T}(z))) >0.$ We are allowed to conclude that for $T$ large enough, for all $z \in {\mathcal D}_3,$ $ \tilde{\check{R}}_{1,T}(z)\in {\mathbb C} \setminus ]- \infty, 0].$ \demend
\ \\
It remains to prove point 4 of proposition \ref{energie-prop-holo}. This is given by the following lemma.

\begin{lem}\label{energie-encadremet-r}
Denote
$$
\cD_4 = \left\{ z\in {\mathbb C}~,~ Re(z) > -\delta_H\theta(1+ \varepsilon) \right\}.
$$
Then, there exist a constant $C$  depending only on $\theta$ and the index $H,$ such that for $T $ large enough,   and $z \in \cD_4$,
\begin{align*}
\sqrt{ \frac{1}{4} }
\leq \left|\exp\Bigl( -\tilde{\check{R}}_T(z) \Bigr) \right|\leq  \sqrt{C \left( 1 + \frac{ 1}{T} + \frac{ 1}{T \varepsilon} \right)}.
\end{align*}
\end{lem}

\noindent{\bf Proof of lemma \ref{energie-encadremet-r}~: } Observe that
\begin{align*}
\left| \exp \Bigl( - \tilde{\check{ R}}_T(z)\Bigr) \right| = \sqrt{ \frac{|\tilde{\check{R}}_{1,T}(z)|}{ (2+p_H) |z| | z+ \delta_H \theta|}}.
\end{align*}
From (\ref{E-DEFR1CT}) and (\ref{E-DEFP}),
\begin{align*}
 \frac{ |\tilde{P}_T(z)|}{(2+p_H)|z| | z + \delta_H \theta |} - \frac{C}{T^2} \leq 
 \frac{|\tilde{\check{R}}_{1,T}(z)|}{ (2+p_H) |z| | z + \delta_H\theta|} \leq 
 \frac{ |\tilde{P}_T(z)|}{(2+p_H)|z| | z + \delta_H \theta |} + \frac{C}{T^2}.
\end{align*}
On the one hand, using the very definition of $\tilde{P}_T,$
\begin{align*}
\frac{\tilde{P}(z)}{(2+p_H) z( z+ \delta_H \theta)}= \frac{1}{z} \left[ z + \frac{r_1}{T} + \frac{\theta r_1(1 - \delta_H)}{T(z+ \delta \theta)}\right]
\end{align*}
and since $r_1 \theta >0,$ then for $z \in {\mathbb C},~~ Re(z) > -\delta_H \theta,$ for $T > 2r_1\theta^{-1} \delta_H^{-1}$
\begin{align*}
\left|\frac{\tilde{P}(z)}{(2+p_H) z( z+ \delta_H \theta)}\right| \geq \frac{1}{|z|} \sqrt{ \left[ Re(z) + \frac{r_1}{T} \right]^2 + Im(z)^2 } \geq \frac{1}{2}.
\end{align*}
On the other hand,
\begin{align*}
\frac{\tilde{P}(z)}{(2+p_H) z( z+ \delta_H \theta)}= 1 + \frac{ r_1}{T} \frac{1}{(2+ p_H) z} + \frac{ r_1 \theta(1 +\delta_H)}{T}\frac{1}{(2+p_H) z( z + \delta_H \theta)}
\end{align*}
and 
\begin{align*}
\left| \frac{\tilde{P}(z)}{(2+p_H) z( z+ \delta_H \theta)} \right| \leq 1 + \frac{ r_1}{\theta \delta_H T(2+ p_H)} +  \frac{ r_1 \theta(1 +\delta_H)}{T}\frac{1}{(2+p_H)  \delta_H^2 \theta^2 \varepsilon )}
\end{align*}
\demend
\ \vspace{1ex}\\
The proof of Proposition \ref{energie-prop-holo}
follows from the conjunction of
Lemmas \ref{lem-energie-fi} to \ref{energie-encadremet-r}.\demend

\noindent{\bf Proof of Proposition \ref{energie-tf-l2}~: }Using the decomposition \eqref{maindeconew} page \pageref{maindeconew}, the third point of Proposition \ref{energie-prop-holo}, the fact that for $a \in \cD_1$,
\begin{align*}
0 ~<~ \inf_{u\in \R} \left\{ \exp( \check{R}_T(a+iu) ) \right\}~\leq~ \sup_{u\in \R} \left\{ \exp(\check{R}_T(a+iu) ) \right\} ~<~ \infty.
\end{align*}
and since $Re(a+ iu) = Re(a),$ we only have to bound 
\begin{align*}
u \mapsto \exp(T ( L( a+ iu)-L(a) )).
\end{align*}
The bound is clearly given by inequality \eqref{E-MAJREF}.\demend


\noindent{\bf Proof of Proposition \ref{energie-prop-maj-r-a-d-t}~:} Recall that, from Lemma \ref{energie-lem-ppte-r}, $
\check{R}_T(a)= \tilde{\check{R}}_T(\varphi(a))$ with $\tilde{\check{R}}_T $ given in \eqref{-energie-def-tilde-r}. Moreover, there exists $0<z_T^-<z_T^+<-\delta_H \theta$ such that
\begin{align*}
\tilde{\check{R}}_{1,T}(z)- (2+p_H)(z-z_T^-)(z-z_T^+) = \CO{2}.
\end{align*}
Since, $0<z_T^-<z_T^+<-\delta_H \theta,$ then for all $z \in ]-\delta_H ,+\infty[$
\begin{align*}
\left|\frac{(2+p_H)(z-z_T^-)(z-z_T^+)}{(2+p_H)z(z-\delta_H \theta)}\right| \leq  1.
\end{align*}
Consequently,
\begin{align*}
\exp( -\check{R}_T(a) )
&\leq \sqrt{ 1 + \frac{\tilde{\check{R}}_{1,T}(\varphi(a))-(2+p_H)(\vp(a)-z_T^-)(\vp(a)-z_T^+)} {(2+p_H)\vp(a)(\vp(a)+\delta_H\theta)}},\\
&= {\cal O} \left( \max \left(1 ; \frac{-1}{T(\varphi(a)+\delta_H \theta)}\right) \right),
\end{align*}
which ends the proof of Proposition \ref{energie-prop-maj-r-a-d-t}.
\demend

\subsection{A contour integral for the Gamma function.}

In order to obtain an asymptotic expansion for $B_T$, it is necessary to make use of the
following lemma which slightly extend Lemma 7.3 of \cite{BR}. First of all, denote by
$f_{a,b}$ the density function of the Gamma $\mathcal{G}(a,b)$ distribution with parameters $a,b>0$, given by
\begin{eqnarray}
\label{fungama}
f_{a,b} (x)= \left \{ \begin{array}{ll}
    {\displaystyle \frac{b^a}{\Gamma(a)} x^{a-1} \exp(- b x)} \ \ \text{ if } \ \ x>0, \vspace{1ex}\\
    0 \ \ \text{ otherwise.}
   \end{array}  \right.
\end{eqnarray}
For all integers $k, \ell \geq 0$, and for all positive real numbers $\sigma^2, \gamma, \nu $, let
$$
v_{k} (a,b,\ell) =\frac{2\pi\sigma ^{2k} i^{\ell}}{2^k k! \gamma^{2k +\ell+1}} f_{a,b}^{(2k+\ell)}(1).
$$
\begin{lem}
\label{leminteg}
For any integers $p>0$ and $\ell \geq 0$, we have
\begin{equation}
\label{valueint}
\int_{\R} \exp\left(-i\gamma u -\frac{\sigma^2u^2}{2T}\right)\frac{u^{\ell}}{(1-2i\nu u )^{a}}\; du
=  \sum_{k=0}^p \frac{v_{k}(a,b,\ell)}{T^k} + {\cal O}\Bigl(\frac{1}{T^{p+1}}\Bigr)
\end{equation}
with $b=\gamma/(2\nu)$.
\end{lem}
{\bf Proof of Lemma \ref{leminteg}.} Denote by $N_{\sigma}$ the Gaussian kernel with
positive variance $\sigma^2$
$$
N_{\sigma}(x)=\frac{1}{\sigma \sqrt{2 \pi}}\exp\left(-\frac{x^2}{2\sigma^2}\right).
$$
It is well-known that
the characteristic functions of $N_{\sigma}$ and $f_{a,b}$ are respectively given by
$$
\widehat{N}_{\sigma}(x)=\exp\left(-\frac{\sigma^2 x^2}{2}\right)
\hspace{1cm}\text{and}\hspace{1cm}
\widehat{f}_{a,b}(x)= \left(1- \frac{ix}{b}\right)^{-a}.
$$
Then, it follows from Rudin \cite{RUD} page 177 that for all integer $\ell \geq 0$ and
for all positive real number $a,b,\tau$
\begin{equation}
\label{profor}
\int_{\R} \exp\left( -ivx - \frac{\tau^2 v^2}{2}\right)
v^{\ell} \widehat{f}_{a,b} (v) \;dv = 2 \pi i^{\ell} f_{a,b} * N_{\tau}^{(\ell)} (x).
\end{equation}
Via the same lines as in \cite{BR} page 17, it is not hard to see that for any $p >0$
\begin{equation}
\label{mainconv}
f_{a,b}*N_{\tau}^{(\ell)} (x)=\sum_{k=0}^p
\frac{\tau^{2k}}{2^k k!}f_{a,b} ^{(2k+\ell)}(x) +  {\cal O}(\tau^{2(p+1)}).
\end{equation}
Hence, we deduce from the conjunction of \eqref{profor} and \eqref{mainconv} with
$\tau^2=\sigma^2/(T\gamma^2)$ that
\begin{equation}
\label{newprofor}
\int_{\R} \exp\left( -ivx - \frac{\sigma^2 v^2}{2T\gamma^2}\right)
v^{\ell} \widehat{f}_{a,b} (v) \;dv = 2 \pi i^{\ell}
\sum_{k=0}^p \frac{\sigma^{2k}}{2^k k!\gamma^{2k}T^k}f_{a,b} ^{(2k+\ell)}(x) +
{\cal O}\Bigl(\frac{1}{T^{p+1}}\Bigr).
\end{equation}
Finally, by taking the values $x=1$ and $b=\gamma/(2\nu)$ together with the
change of variables $u=v/\gamma$ in \eqref{newprofor}, we find that
\begin{equation}
\int_{\R}\exp\left(-i\gamma u  -\frac{\sigma^2u^2}{2T}\right)
\frac{u^{\ell}}{(1- 2i\nu u)^{a}}\; du=
\sum_{k=0}^p \frac{v_{k}(a,b,\ell)}{T^k} +
{\cal O}\Bigl(\frac{1}{T^{p+1}}\Bigr)
\end{equation}
which completes the proof of Lemma \ref{leminteg}.\demend

\bibliographystyle{plain}
\bibliography{Bercu_Coutin_Savy}

\begin{thebibliography}{10}

\bibitem{BaR}
R.~R. Bahadur and R.~Ranga~Rao.
\newblock On deviations of the sample mean.
\newblock {\em Ann. Math. Statist.}, 31:1015--1027, 1960.

\bibitem{BGR}
B.~Bercu, F.~Gamboa, and A.~Rouault.
\newblock Large deviations for quadratic forms of stationary {G}aussian
  processes.
\newblock {\em Stochastic Process. Appl.}, 71(1):75--90, 1997.

\bibitem{BR}
B.~Bercu and A.~Rouault.
\newblock Sharp large deviations for the {O}rnstein-{U}hlenbeck process.
\newblock {\em Theory Probab. Appl.}, 46(1):1--19, 2002.

\bibitem{B08}
J.~P.~N. Bishwal.
\newblock Large deviations in testing fractional {O}rnstein-{U}hlenbeck models.
\newblock {\em Statist. Probab. Lett.}, 78(8):953--962, 2008.

\bibitem{BK}
A.~Brouste and M.~L. Kleptsyna.
\newblock Asymptotic properties of mle for partially observed fractional
  diffusion system.
\newblock Preprint, 2008.

\bibitem{BD}
W.~Bryc and A.~Dembo.
\newblock Large deviations for quadratic functionals of {G}aussian processes.
\newblock {\em J. Theoret. Probab.}, 10(2):307--332, 1997.
\newblock Dedicated to Murray Rosenblatt.

\bibitem{DZ}
A.~Dembo and O.~Zeitouni.
\newblock {\em Large deviations techniques and applications}, volume~38 of {\em
  Applications of Mathematics (New York)}.
\newblock Springer-Verlag, New York, second edition, 1998.

\bibitem{FlP}
D.~Florens-Landais and H.~Pham.
\newblock Large deviations in estimation of an {O}rnstein-{U}hlenbeck model.
\newblock {\em J. Appl. Probab.}, 36(1):60--77, 1999.

\bibitem{GRZ}
F.~Gamboa, A.~Rouault, and M.~Zani.
\newblock A functional large deviations principle for quadratic forms of
  {G}aussian stationary processes.
\newblock {\em Statist. Probab. Lett.}, 43(3):299--308, 1999.

\bibitem{KLB1}
M.~L. Kleptsyna and A.~Le~Breton.
\newblock Statistical analysis of the fractional {O}rnstein-{U}hlenbeck type
  process.
\newblock {\em Stat. Inference Stoch. Process.}, 5(3):229--248, 2002.

\bibitem{Leb}
N.~N. Lebedev.
\newblock {\em Special functions and their applications}.
\newblock Revised English edition. Translated and edited by Richard A.
  Silverman. Prentice-Hall Inc., Englewood Cliffs, N.J., 1965.

\bibitem{NVV}
I.~Norros, E.~Valkeila, and J.~Virtamo.
\newblock An elementary approach to a {G}irsanov formula and other analytical
  results on fractional {B}rownian motions.
\newblock {\em Bernoulli}, 5(4):571--587, 1999.

\bibitem{PKS1}
B.~L.~S. Prakasa~Rao.
\newblock Sequential estimation for fractional ornstein-uhlenbeck type process.
\newblock {\em Sequential Anal.}, 23(1):33--44, 2004.

\bibitem{PKS2}
B.~L.~S. Prakasa~Rao.
\newblock Estimation for translation of a process driven by fractional brownian
  motion.
\newblock {\em Stoch. Anal. Appl.}, 23(6):1199--1212, 2005.

\bibitem{RUD}
W.~Rudin.
\newblock {\em Real and complex analysis}.
\newblock McGraw-Hill Book Co., New York, third edition, 1987.

\end{thebibliography}

\end{document}